\documentclass[11pt]{amsart}
\usepackage[T1]{fontenc}
\usepackage{lmodern,amsmath,amssymb,amsthm,mathtools}
\usepackage[margin=1.08in]{geometry}
\usepackage{microtype,booktabs,array,xcolor}
\usepackage[colorlinks=true,linkcolor=blue!50!black,citecolor=blue!50!black,urlcolor=blue!50!black]{hyperref}
\hypersetup{pdftitle={Rational Elliptic Period Relations Across Realizations:
Complex Lattices, Local p-adic Periods, and Drinfeld Modules},pdfauthor={Mohammadreza Mohajer},hypertexnames=false}
\newtheorem{theorem}{Theorem}[section]
\newtheorem{proposition}[theorem]{Proposition}
\newtheorem{lemma}[theorem]{Lemma}
\newtheorem{corollary}[theorem]{Corollary}
\theoremstyle{definition}

\newtheorem{example}[theorem]{Example}
\newtheorem{question}[theorem]{Question}
\newtheorem{conjecture}[theorem]{Conjecture}
\theoremstyle{remark}
\newtheorem{remark}[theorem]{Remark}
\numberwithin{equation}{section}
\newcommand{\Q}{\mathbb Q}
\newcommand{\Z}{\mathbb Z}
\newcommand{\C}{\mathbb C}
\newcommand{\F}{\mathbb F}
\newcommand{\Pj}{\mathbb P}
\newcommand{\Qp}{\Q_p}
\newcommand{\Qpb}{\overline{\Q}_p}
\newcommand{\Cp}{\C_p}
\newcommand{\BdR}{B_{\mathrm{dR}}}
\newcommand{\BdRp}{B_{\mathrm{dR}}^+}
\newcommand{\HdR}{H^1_{\mathrm{dR}}}
\newcommand{\DdR}{D_{\mathrm{dR}}}
\newcommand{\Dcris}{D_{\mathrm{cris}}}
\newcommand{\trdeg}{\operatorname{trdeg}}
\newcommand{\End}{\operatorname{End}}
\newcommand{\Lie}{\operatorname{Lie}}
\newcommand{\GL}{\operatorname{GL}}
\newcommand{\Fil}{\operatorname{Fil}}
\newcommand{\Gm}{\mathbb G_m}
\newcommand{\Ga}{\mathbb G_a}
\newcommand{\HH}{\mathfrak H}
\DeclareMathOperator{\Ann}{Ann}

\DeclareMathOperator{\Hom}{Hom}

\DeclareMathOperator{\MT}{MT}

\newcommand{\Qb}{\overline{\Q}}
\newcommand{\Bcris}{B_{\mathrm{cris}}}
\title[Rational Elliptic Period Relations Across Realizations]{Rational Elliptic Period Relations Across Realizations:\\
Complex Lattices, Local $p$-adic Periods, and Drinfeld Modules}
\author{Mohammadreza Mohajer}
\address{Department of Mathematics, Physics and Geology, Cape Breton University,
1250 Grand Lake Road, Sydney, NS B1P 6L2, Canada}
\email{mohammadreza.mohajer96@gmail.com}
\email{Mohammadreza\_Mohajer@cbu.ca}
\begin{document}
\begin{abstract}
We determine the algebraic relations among full de Rham comparison
periods of elliptic curves over finite extensions of $\Qp$. The
connected local monodromy group gives the transcendence degree and
an explicit affine linear relation ideal for the full matrix in
Fontaine's field $\BdR$. After normalizing the holomorphic row and
reducing to $\Cp$, the transcendence degree drops by one and the
normalized determinant supplies the additional relation. Every
polynomial relation of degree at most $d$ is detected modulo
$t^{d+1}$. For $d\ge1$ this bound is sharp for each elliptic curve.

A rational quasi-period row gives a Galois-stable Tate line.
Together with an algebraic holomorphic period ratio it forces
potential splitting, which identifies the exact ordinary CM
obstruction. We compute the full Tate matrix. For ordinary
non-CM curves with $p\ge5$, the standard quasi-period row is rational precisely
at a zero of Katz's $E_{2,p}$. Away from the zeros the two full
ratios are algebraically independent. At the zeros, a nonzero
quasi-period and the normalized holomorphic ratio have algebraically
independent images in $\Cp$. The complex question that motivates
these results, the analogue for Drinfeld modules, and a conditional
arithmetic consequence for Katz values complete the comparison.
\end{abstract}
\maketitle
\section{Introduction}\label{sec:intro}
\subsection{Question A and its period-theoretic motivation}
Let $\Lambda=\Z\omega_1+\Z\omega_2\subset\C$ be a complex
lattice, with the basis chosen so that
$\tau=\omega_2/\omega_1\in\HH$, where
$\HH=\{z\in\C:\operatorname{Im}z>0\}$. Let $\wp_\Lambda$
denote the Weierstrass elliptic function associated with $\Lambda$.
The corresponding Weierstrass zeta function $\zeta_\Lambda$ satisfies
$\zeta_\Lambda'=-\wp_\Lambda$ and is normalized by
$\zeta_\Lambda(z)=z^{-1}+O(z^3)$ as $z\to0$.
Its additive increments are the quasi-periods:
\[
 \zeta_\Lambda(z+\omega_i)-\zeta_\Lambda(z)=\eta_i,\qquad
 \eta_1\omega_2-\eta_2\omega_1=2\pi i.
\]
Let $g_2=g_2(\Lambda)$ and $g_3=g_3(\Lambda)$ be the
Weierstrass invariants, and put $F=\Q(g_2,g_3)$. The elliptic curve
associated with $\Lambda$ has the short Weierstrass model
\begin{equation}
 E_\Lambda:\quad y^2=4x^3-g_2x-g_3,\qquad
 \omega=\frac{dx}{y},\qquad \xi=-x\frac{dx}{y}.
 \label{eq:complex-model}
\end{equation}
The classes of $\omega$ and the second-kind differential $\xi$
form an $F$-basis of $H^1_{\mathrm{dR}}(E_\Lambda/F)$.
Under the Weierstrass uniformization, these differentials pull back
to $dz$ and $d\zeta_\Lambda(z)$, respectively.
Let $\gamma_1,\gamma_2$ be the integral homology basis of
$E_\Lambda(\C)$ corresponding to $\omega_1,\omega_2$.
Integration gives the de Rham--Betti comparison isomorphism
\[
 H^1_{\mathrm{dR}}(E_\Lambda/F)\otimes_F\C
 \xrightarrow{\sim}
 H^1_{\mathrm B}(E_\Lambda(\C),\Q)\otimes_\Q\C.
\]
We arrange its periods with de Rham classes indexing rows and
homology cycles indexing columns:
\begin{equation}
 P_\Lambda=
 \begin{pmatrix}
  \displaystyle\int_{\gamma_1}\omega&
  \displaystyle\int_{\gamma_2}\omega\\[4pt]
  \displaystyle\int_{\gamma_1}\xi&
  \displaystyle\int_{\gamma_2}\xi
 \end{pmatrix}
 =
 \begin{pmatrix}
  \omega_1&\omega_2\\
  \eta_1&\eta_2
 \end{pmatrix}.
 \label{eq:complex-matrix}
\end{equation}
The first row is the \emph{holomorphic period row}, and the
second is the \emph{quasi-period row}.

We call the quasi-period row \emph{rational} when its projective
class satisfies
\[
[\eta_1:\eta_2]\in\Pj^1(\Q)\subseteq\Pj^1(\C).
\]
Thus $[\eta_1:\eta_2]\in\Pj^1(\Q)$ means that
$(\eta_1,\eta_2)=c(a,b)$ for some $c\in\C^\times$ and
$(a,b)\in\Q^2\setminus\{(0,0)\}$. This includes a zero quasi-period.
The Legendre relation excludes simultaneous vanishing.

The central question is
\begin{equation}
 [\eta_1:\eta_2]\in\Pj^1(\Q)
 \quad\Longrightarrow\quad \tau\notin\Qb.
 \tag{A}\label{eq:A}
\end{equation}

We call~\eqref{eq:A} Question A. Kumar posed the necessity question
in~\cite[\S4]{Kumar} and gave the modular reformulation recalled in
Proposition~\ref{prop:kumar}.

Question A is motivated by the expectation that algebraic relations
among periods arise from the geometry of the underlying elliptic curve. Grothendieck's
period conjecture \cite[p.~102, footnote~10]{Grothendieck} expresses this principle for varieties over
algebraic numbers; Andr\'e's generalized formulation includes the
field of definition when the coefficients are transcendental
\cite{Andre,Bertolin20}. In its elliptic form the generalized
conjecture predicts
\begin{equation}
 \trdeg_\Q F(\omega_1,\omega_2,\eta_1,\eta_2)
 \ \geq\ \dim\MT(E)
 =
 \begin{cases}4,&E\text{ without CM},\\2,&E\text{ with CM},\end{cases}
 \label{eq:motivic-intro}
\end{equation}
for an elliptic curve $E$ defined over a finitely generated subfield $F\subset\C$
and a de Rham basis defined over $F$. Here $\MT(E)$ is the
Mumford--Tate group of the rational Hodge structure $H^1(E(\C),\Q)$.
The transcendence degree is taken over $\Q$, so it includes the
contribution from $F$. When $F$ is algebraic over $\Q$, the classical
period conjecture predicts equality; for general $F$, the generalized
formulation gives the lower bound in~\eqref{eq:motivic-intro}.
Proposition~\ref{prop:motivic-A} proves that this generalized bound
implies Question A. This application requires the generalized formulation:
Proposition~\ref{prop:zero-j} shows that every lattice satisfying
$[\eta_1:\eta_2]\in\Pj^1(\Q)$ has transcendental $j$-invariant.
The associated elliptic curve therefore has no model over $\Qb$,
so the classical period conjecture for curves over algebraic
numbers does not apply directly.

Kumar's reformulation, Proposition~\ref{prop:kumar}, turns Question A
into the transcendence problem for the zeros of the holomorphic
weight-two Eisenstein series $E_2$. We examine that implication in other
realizations, determine which parts admit unconditional proofs,
and identify arithmetic information relevant to the complex question.

\subsection{The realizations and the results}\label{subsec:results}
A period relation depends on a cycle space, a de Rham class, and
a comparison pairing. Over a finite extension $K/\Qp$, the cycle
space is the rational Tate module $V_p(E)$ and the pairing takes
values in Fontaine's de Rham period field $\BdR$. We call these
periods \emph{full} because we take the entries of the comparison
matrix in $\BdR$ before reducing by its filtration. Write $t$ for
the cyclotomic period, a uniformizer of $\BdRp$.

\begin{enumerate}
\item \emph{Rational rows and local splitting.}
Lemma~\ref{lem:rows} identifies a rational period row with a
Galois-stable Tate line. Lemma~\ref{lem:splitting} shows that a
rational quasi-period row and an algebraic holomorphic period
point together force the Tate representation to split over a
finite extension. Theorem~\ref{thm:local-A} proves the local
implication for curves without geometric CM.
Corollary~\ref{cor:exact-obstruction} identifies the precise CM
exception, and Proposition~\ref{prop:Betti-CM} describes its
dependence on the rational structure on cycles.

\item \emph{Full and reduced period fields.}
Theorem~\ref{thm:dimensions} gives the five possibilities for the
connected local monodromy and the transcendence degree of the
full period matrix. Corollary~\ref{cor:local-relation-ideal} gives
the corresponding affine linear relation ideals.
Theorem~\ref{thm:graded-relation-ideal} determines the exact
relations of the normalized matrix in $\Cp$: the transcendence
degree drops by one, and in the indicated coordinates the
normalized determinant gives the additional equation.
Corollary~\ref{cor:bounded-degree-truncation} shows that direct
evaluation in $B_{d+1}=\BdRp/t^{d+1}\BdRp$ detects every
polynomial relation of degree at most $d$ with coefficients
algebraic over $\Qp$. For $d\ge1$ the bound is sharp for every elliptic curve. Proposition~\ref{prop:holomorphic-graded} identifies the
algebraic linear relations of the full holomorphic row with
those of its degree-one reduction.

\item \emph{Tate curves and ordinary reduction.}
Theorem~\ref{thm:tate-matrix} computes the full Tate matrix in the
standard Weierstrass basis. Theorem~\ref{thm:tate-independence}
proves independence of its two ratios, and
Corollary~\ref{cor:correction} identifies the unique correction
of the second-kind class with a rational row. The same Kummer
argument gives the rank formula for toric $1$-motives in
Theorem~\ref{thm:Kummer-rank}.
For an ordinary elliptic curve, Theorem~\ref{thm:ordinary-condition}
identifies rationality of the standard quasi-period row with
vanishing of Katz's $E_{2,p}$. Corollary~\ref{cor:ordinary-ratios}
gives algebraic independence of the two full ratios away from the
zeros on the non-CM ordinary locus.
Corollary~\ref{cor:ordinary-zero-periods} determines the full and
reduced period fields in a unit-root basis. At a Katz zero the
standard class is already a unit-root class. The example over
$\Q_{11}$ applies these results to a zero supplied by the
classification of Sebbar--Besrour \cite{SBp}.
\end{enumerate}

For the equal-characteristic comparison, put $A=\F_q[\theta]$ and
$k=\F_q(\theta)$, where $q$ is a power of $p$. Rank-two Drinfeld
modules of generic $A$-characteristic have lattices and entire
quasi-periodic functions. Section~\ref{sec:drinfeld} gives the
zero reformulation for these lattices and applies Yu and
Chang--Papanikolas over the algebraic closure of $k$.
The unrestricted complex and Drinfeld questions remain open.
Proposition~\ref{prop:transfer} specifies the comparison
compatibility that would transport the two marked complex
relations. Under the ordinary case of the crystalline period
conjecture of Ancona--Fr\u{a}\c{t}il\u{a}~\cite[Conjecture~8.3]{AF},
Proposition~\ref{prop:ordinary-arithmetic} proves that the two
Frobenius eigenline coordinates, one of which is $E_{2,p}/12$,
are algebraically independent over $\Qb$. In particular this
conjecture implies that $j_{11}$ is transcendental over $\Q$.

\subsection{Related work}\label{subsec:related}
Serre \cite[Chapter~IV, Appendix~A]{Serre} determines the local
monodromy groups used in Theorem~\ref{thm:dimensions}. His
Hodge--Tate orbit formula
\cite[\S1.5, Theorem~3]{SerreHT} calculates the transcendence
degree of the grading through the centralizer of its cocharacter.
Howe \cite[Theorem~1.2]{Howe} studies the included Hodge--Tate
filtration line, and Howe--Klevdal
\cite[Remark~1.1.4]{HoweKlevdal} distinguish that line from the
canonical grading for ordinary elliptic curves. In this case
algebraicity of the grading is equivalent to torsion of the
Serre--Tate parameter. Our degree-one row is the covector giving
projection onto the weight-one summand, as explained in
Subsection~\ref{subsec:holomorphic-graded}. We use the same
monodromy groups to describe the full comparison matrix and
both rows of its normalized reduction.

Berger \cite[\S II.4]{Berger} describes the Tate representation,
its Kummer class, and the corresponding filtered modules.
Theorem~\ref{thm:tate-matrix} identifies the periods of the
specific Weierstrass classes using a theta primitive and the
geometric comparison for semistable abelian varieties
\cite{ColemanIovita}. For ordinary curves, Katz
\cite[\S5.7]{Katz} identifies the unit-root coefficient
$E_{2,p}(E,\omega)/12$, and Mazur--Stein--Tate \cite{MST}
compute it through Frobenius. The same coefficient appears in
the canonical $p$-adic sigma function of Mazur--Tate
\cite{MT}. Norman \cite{Norman} gives a theta construction,
Blakestad--Grant \cite{BG} develop universal sigma and zeta
functions, and Balakrishnan--Besser \cite{BB} relate the sigma
function to local heights.

Sebbar--Besrour \cite[Sections~3 and~8]{SBp} count the ordinary
Katz zeros, determine their unramified fields of definition,
classify the CM zeros, and study isogenies among the zeros.
Section~\ref{sec:ordinary} applies this classification to full
comparison periods. Their Remark~8.11 and Proposition~8.12 give
the source for the conditional nonvanishing argument revisited in
Proposition~\ref{prop:ordinary-arithmetic}. The period space of
Ancona--Fr\u{a}\c{t}il\u{a} \cite{AF} concerns endomorphisms of
the reduction expressed in a global algebraic de Rham basis.
Subsection~\ref{subsec:crystalline-consequence} explains why the
base field there is $\Qb$, whereas the full local matrix is
studied over $\Qpb$.

Kumar \cite{Kumar} relates rational quasi-period dependence to
the modular orbits of $E_2$-zeros. Sebbar
\cite[Sections~3--4 and~13--14]{Sebbar} studies the zeros,
transcendence of their invariants, and conditional period
arguments. The universal vector extension provides an
algebraic-group interpretation of the complex relation.
W\"ustholz \cite{Wustholz} and Huber--W\"ustholz \cite{HW}
describe linear relations among periods of $1$-motives.
Fuchs--Pham \cite{FP15,FP16} and Pham \cite{Pham}
treat $p$-adic logarithms and analytic Weierstrass zeta values
of algebraic points. Sections~\ref{sec:complex}
and~\ref{sec:analytic} apply the subgroup theorems to their
respective period and logarithm arguments.
The generalized period conjectures for $1$-motives and their
formulation through elliptic integrals and semi-elliptic
functions are developed in
\cite{Bertolin02,Bertolin20,BW,Bertolin25}.

Mohajer--Sebbar \cite{MS} construct a $B_2$-valued integration
pairing for good-reduction $1$-motives and specify the
cycle and de Rham arguments for their number-field period
spaces. Proposition~\ref{prop:MS-intersection} compares their
cycle arguments with the constant Tate vectors of the present
paper. Our filtration proof treats the full comparison and
applies to all elliptic reduction types.

In generic $A$-characteristic, Gekeler's de Rham comparison
\cite{Gekeler} gives the matrix description and his false
Eisenstein formula \cite[Theorem~7.10]{GekelerQuasi} gives
the zero condition. Yu's quasi-period theorem \cite{Yu90}
and sub-$t$-module theorem \cite{Yu97} give the setting for
linear relations. Chang--Papanikolas \cite{CP} determine the
transcendence degree of the matrix, and Chang \cite{ChangQuasi}
studies nonzero quasi-modular values at algebraic moduli.
Corollary~\ref{cor:drinfeld-j} treats the $j$-invariant at a
zero of the standard quasi-period function.

\subsection{Conventions and organization}
In the complex sections, algebraicity means algebraicity over
$\Q$. In the local sections it means algebraicity over $\Qp$,
unless the base $\Q$ is explicit. Thus an element of $K$ is
algebraic over $\Qp$ and may be transcendental over $\Q$.
We write $\Qpb$ for the chosen algebraic closure of $\Qp$,
embedded in $\BdR$ as in Section~\ref{sec:local}, and $\Cp$ for
its completion with respect to the $p$-adic absolute value.
Geometric CM means $\End_{\overline K}(E)\ne\Z$.
Formal CM concerns the formal $p$-divisible group and is defined
before Table~\ref{tab:local}. Rationality and algebraicity of a
period row refer to its projective class. Unless another field is
specified, the reference field is $\Q$ in the complex setting,
$\Qp$ in the local setting, and $k=\F_q(\theta)$ in the Drinfeld
setting.

Sections~\ref{sec:complex}--\ref{sec:monodromy} establish the
complex reformulation and the local structural results.
Sections~\ref{sec:tate}--\ref{sec:ordinary} give explicit
calculations. Section~\ref{sec:drinfeld} explains the
equal-characteristic comparison, and Section~\ref{sec:analytic}
connects analytic primitives with de Rham classes.
Section~\ref{sec:bridge} returns to the motivic motivation
and the unresolved arithmetic questions.

\subsection*{Acknowledgments}
I am grateful to Abdellah Sebbar for introducing me to Question~\eqref{eq:A}
and to the connection between elliptic quasi-periods and Eisenstein
series. Our discussions prompted me to study analogous period relations
over $p$-adic local fields and for Drinfeld modules.

\section{The complex formulation}\label{sec:complex}
\subsection{Normalization and the zero problem}
We use the short Weierstrass model and de Rham basis
$\omega=dx/y$, $\xi=-x\,dx/y$ introduced
in~\eqref{eq:complex-model}. Their period matrix is
\eqref{eq:complex-matrix}.
For $\Lambda=\omega_1(\Z+\Z\tau)$,
\begin{equation}
 \eta_1=\frac{\pi^2}{3\omega_1}E_2(\tau),\quad
 g_2=\frac{4\pi^4}{3\omega_1^4}E_4(\tau),\quad
 g_3=\frac{8\pi^6}{27\omega_1^6}E_6(\tau).
 \label{eq:complex-normalizations}
\end{equation}
Here $q=e^{2\pi i\tau}$ and $\sigma_m(n)=\sum_{d\mid n}d^m$.
We use
$E_2(q)=1-24\sum_{n\geq1}\sigma_1(n)q^n$,
$E_4(q)=1+240\sum_{n\geq1}\sigma_3(n)q^n$, and
$E_6(q)=1-504\sum_{n\geq1}\sigma_5(n)q^n$.
The second-kind differential $\xi$ has a double pole and zero
residue at the origin. Its de Rham class is well defined,
and $\omega,\xi$ form a cohomology basis.

The following reformulation is due to Kumar~\cite{Kumar}; see
\cite[Propositions~3.1 and~3.3]{Sebbar} for the formulation in terms
of rational projective rows and changes of lattice basis.
\begin{proposition}[Kumar's reformulation]\label{prop:kumar}
Question~\eqref{eq:A} is equivalent to the transcendence over $\Q$
of every zero of $E_2$ in $\HH$.
\end{proposition}
\begin{proof}
A rational dependence gives coprime integers $a,b$, not both zero,
with $a\eta_1+b\eta_2=0$. The period
$\omega'_1=a\omega_1+b\omega_2$ is primitive and extends to a
positively oriented integral basis of $\Lambda$. After scaling $\omega'_1$ to $1$,
equation~\eqref{eq:complex-normalizations} gives $E_2(\tau')=0$.
The new ratio $\tau'$ is an integral fractional linear transform
of $\tau$, and this transformation preserves algebraicity in both
directions. Conversely, a zero gives a lattice with $\eta'_1=0$.
This is the argument of \cite{Kumar} and
\cite[Propositions~3.1 and~3.3]{Sebbar}.
\end{proof}

We recall the following consequence of Nesterenko's theorem,
proved in \cite[Proposition~4.3]{Sebbar}.
\begin{proposition}[The coefficient field at a zero]\label{prop:zero-j}
If $E_2(\rho)=0$ in $\HH$, then $j(\rho)$ is transcendental over
$\Q$. A lattice with rationally dependent quasi-periods
therefore has transcendental $j$ and has no CM.
\end{proposition}
\begin{proof}
By \cite[Proposition~4.3]{Sebbar}, Nesterenko's theorem
\cite{Nesterenko} makes $q,E_4(\rho),E_6(\rho)$ algebraically
independent when $E_2(\rho)=0$. The formula
$j=1728E_4^3/(E_4^3-E_6^2)$ then excludes an algebraic $j$.
Proposition~\ref{prop:kumar} transfers this conclusion to a lattice
with rationally dependent quasi-periods. CM $j$-invariants are
algebraic, so such a lattice has no CM.
\end{proof}
This is why the period conjecture only for curves over $\Qb$
does not directly settle Question A.

\subsection{Analytic subgroups and algebraic period columns}
\label{subsec:complex-subgroups}
The universal vector extension is the algebraic group on which
we apply the complex analytic subgroup theorem. Its $p$-adic
logarithm is used in Subsection~\ref{subsec:analytic-subgroups}.

In the form of W\"ustholz's theorem used here
\cite{Wustholz}\cite[Theorem~6.2]{HW}, let $G/\Qb$ be a connected
commutative algebraic group and let $L\subseteq\Lie(G)$ be a
$\Qb$-linear subspace. If
\[
 0\ne u\in L\otimes_{\Qb}\C,\qquad \exp_G(u)\in G(\Qb),
\]
then a connected algebraic subgroup $H\subseteq G$ over $\Qb$
satisfies
\[
 0\ne\Lie(H)\subseteq L,\qquad u\in\Lie(H)\otimes_{\Qb}\C.
\]
The theorem applies to period vectors, whose exponential is the
identity.

For an elliptic curve $E$ over an algebraically closed field of
characteristic zero, identify the one-dimensional vector kernel
of its universal vector extension with $\Ga$ and write
\begin{equation}
 0\longrightarrow\Ga\longrightarrow E^\natural
   \xrightarrow{\pi}E\longrightarrow0.
 \label{eq:universal-extension}
\end{equation}
This is the extension whose nonzero class generates
$\operatorname{Ext}^1(E,\Ga)=H^1(E,\mathcal O_E)$; every extension
of $E$ by a vector group is obtained from it by pushout.

\begin{lemma}[Subgroups of the universal vector extension]
\label{lem:universal-subgroups}
Every proper connected algebraic subgroup of $E^\natural$ is
contained in the additive kernel of $\pi$.
\end{lemma}
\begin{proof}
A positive-dimensional proper connected subgroup $H$ has dimension
one. If $\pi(H)\ne0$, then $\pi|_H$ is surjective with finite
kernel. A connected affine algebraic group has no nonconstant
group homomorphism to an elliptic curve, so the classification of
one-dimensional connected groups makes $H$ an elliptic curve.
Thus $f=\pi|_H:H\to E$ is an isogeny. The inclusion of $H$ into
$E^\natural$ splits the pullback of \eqref{eq:universal-extension}
along $f$. Pullback on $H^1(\mathcal O)$ is injective in
characteristic zero, since its composition with trace is
multiplication by $\deg f$. It cannot kill the nonzero extension
class. Hence $\pi(H)=0$. The zero-dimensional connected subgroup
is trivial and satisfies the same conclusion.
\end{proof}

Suppose now that $E/\Qb$, and choose a de Rham basis
$(\omega,\xi)$ over $\Qb$, with $\omega$ holomorphic.
Using the canonical principal polarization, the de Rham
realization identifies $\Lie(E^\natural)$ with
$H^1_{\mathrm{dR}}(E/\Qb)^*$. Evaluation against $\omega$
identifies $\Lie(E)$ with $\Qb$. Under these identifications,
$d\pi$ is restriction to the holomorphic line $\Qb\omega$,
hence the first coordinate. The complex exponential lattice
consists of the vectors
\[
 u_\gamma=\left(\int_\gamma\omega,\int_\gamma\xi\right),
 \qquad \gamma\in H_1(E(\C),\Z).
\]
See \cite[Chapters~4--5]{HW} for this realization of the universal
vector extension.

\begin{proposition}[Independence within a complex period column]
\label{prop:complex-column}
For every nonzero $\gamma\in H_1(E(\C),\Z)$, the two coordinates
of $u_\gamma$ are linearly independent over $\Qb$.
\end{proposition}
\begin{proof}
Write $u_\gamma=(\Omega,H)$. The holomorphic period $\Omega$ is
nonzero, since integration of $\omega$ identifies integral
homology with the elliptic lattice. A relation
$a\Omega+bH=0$, with $a,b\in\Qb$ not both zero, has $b\ne0$.
The algebraic line $L=\ker(aX+bY)$ therefore projects nontrivially
to $\Lie(E)$.

The classical analytic subgroup theorem
\cite[Theorem~6.2]{HW}, applied to $u_\gamma\ne0$ with
$\exp_{E^\natural}(u_\gamma)=0$, gives a connected algebraic
subgroup $H\subseteq E^\natural$ with
$u_\gamma\in\Lie(H)\otimes_{\Qb}\C$ and
$0\ne\Lie(H)\subseteq L$. Since $L$ is a line,
$\Lie(H)=L$. Its nonzero projection to $\Lie(E)$ contradicts
Lemma~\ref{lem:universal-subgroups}.
\end{proof}

\begin{remark}[An analytic-subgroup proof of the coefficient-field obstruction]
\label{rem:complex-j-subgroup}
Proposition~\ref{prop:complex-column} gives a second proof of
Proposition~\ref{prop:zero-j}. If a lattice has algebraic $j$,
a homothety makes both Weierstrass invariants algebraic: when
$g_2\ne0$, normalize $g_2$ to $1$ and use the formula for $j$
to see that $g_3^2$ is algebraic; when $g_2=0$, normalize
$g_3$ to $1$. The standard quasi-periods are multiplied by a
common nonzero factor, so rational dependence is preserved.
Clearing denominators then gives a nonzero integral cycle with
zero $\xi$-period, contrary to
Proposition~\ref{prop:complex-column}.

The Nesterenko proof above gives the obstruction directly in
modular coordinates. The present proof identifies its algebraic
subgroup mechanism. Both proofs show why a hypothetical
counterexample to Question A lies beyond the number-field
hypothesis of the classical analytic subgroup theorem.
\end{remark}

\subsection{The CM qualification}
Having four transcendental entries is weaker than having transcendence
degree four. For an elliptic curve over $\overline\Q$ with a de Rham basis
defined over $\overline\Q$, the expected transcendence degree of the complex period
matrix is $4$ in the non-CM case and $2$ in the CM case. The CM value
is known from Chudnovsky's elliptic period theorems~\cite{Chudnovsky}. A formulation using
only $j$ must still specify the normalization of the de Rham basis.

For the square lattice $\Z+i\Z$ one has
\[
 \omega_1=1,\quad\omega_2=i,\quad
 \eta_1=\pi,\quad\eta_2=-i\pi,\quad g_3=0.
\]
The field generated by $g_2,g_3$ and the four period entries is
contained in $\Q(i,g_2,\pi)$, of transcendence degree at most $2$,
in agreement with the CM bound in~\eqref{eq:motivic-intro}.
Kumar's concluding remarks \cite[\S4]{Kumar} cite the square
lattice as a counterexample to the necessity of transcendence of
$\tau$. In the normalization above the quasi-period ratio is
$-i$, so the row is algebraic and lies outside $\Pj^1(\Q)$.
Thus the example does not satisfy the premise of Question~\eqref{eq:A}.

\section{Full local periods and a splitting argument}\label{sec:local}
Fix a finite extension $K/\Qp$ and an algebraic closure
$\overline K=\Qpb$, with completion $\Cp$. Write
$G_K=\operatorname{Gal}(\overline K/K)$.
Fontaine's ring $\BdRp$ is a complete discrete valuation ring
with residue field $\Cp$; its residue map is
$\theta:\BdRp\to\Cp$. Its fraction field is $\BdR$.
The algebraic field $\Qpb$ embeds in $\BdRp$, and $\theta$
restricts to the given embedding in $\Cp$.
Choose a compatible system of primitive $p^n$-th roots of unity
$\varepsilon=(1,\zeta_p,\zeta_{p^2},\ldots)$ and put
$t=\log[\varepsilon]$. Here brackets denote the multiplicative
Teichm\"uller lift in Fontaine's construction. The element $t$
generates the maximal ideal of $\BdRp$ and satisfies
$\sigma(t)=\chi(\sigma)t$, where
$\chi:G_K\to\Z_p^\times$ is the cyclotomic character.

Let $T_p(E)=\varprojlim_n E[p^n](\overline K)$, with transition
maps multiplication by $p$. Its rationalization
$V=V_p(E)=T_p(E)\otimes_{\Z_p}\Qp$ has dimension two.
For a de Rham representation $R$, our convention is
$\DdR(R)=(\BdR\otimes_{\Qp}R)^{G_K}$.
In particular $\HdR(E/K)=\DdR(V^*)$. Its Hodge filtration is
\[
 \Fil^0\HdR(E/K)=\HdR(E/K),\quad
 \Fil^1\HdR(E/K)=K\omega,\quad
 \Fil^2\HdR(E/K)=0.
\]

The de Rham comparison theorem gives a pairing
\begin{equation}
 V\times\HdR(E/K)\longrightarrow\BdR
 \label{eq:pairing}
\end{equation}
whose extension to $\BdR$ is perfect. It evaluates a cohomology
class, viewed through $\HdR(E/K)=\DdR(V^*)$, on a Tate-module
vector. The pairing is $\Qp$-linear on
$V$, $K$-linear on cohomology, and $G_K$-equivariant.
We use $\BdR^{G_L}=L$ for finite $L/K$, exactness of the de Rham
functor on de Rham representations, and compatibility with cup products
and filtrations; see~\cite{BC,Berger}.

For bases $e_1,e_2$ of $V$ and $\omega,\xi$ of $\HdR(E/K)$,
with $\omega$ holomorphic, set
\begin{equation}
 W_i=\langle e_i,\omega\rangle,\quad
 H_i=\langle e_i,\xi\rangle,\quad
 P=\begin{pmatrix}W_1&W_2\\H_1&H_2\end{pmatrix}.
 \label{eq:matrix}
\end{equation}
This matrix is invertible over $\BdR$. Thus neither period row vanishes.
Filtration compatibility puts $W_i$ in $t\BdRp$ and $H_i$ in
$\BdRp$. A projective point $[W_1:W_2]$ is formed in the field
$\BdR$, even when an entry is zero.
For the short Weierstrass basis
$\omega=dx/y$, $\xi=-x\,dx/y$, choose a symplectic Tate-module
basis whose Weil pairing agrees with the generator of $\Z_p(1)$
specified by $\varepsilon$. We use the Weil-pairing convention
in which $e_n(\zeta,q^{1/n})=\zeta$ on a Tate curve
\cite[\S1.2.7]{NekovarParity}.
In the local parameter $z$ with $\omega=dz$ at the identity,
$\xi=(-z^{-2}+O(z^2))\,dz$. The residue of $z\xi$ is $-1$,
which is the de Rham cup product of $\omega,\xi$ in our
convention. Compatibility with the Weil pairing gives
\begin{equation}
 \det P=-t. \label{eq:legendre-p}
\end{equation}
In arbitrary de Rham and cycle bases $\det P$ is a nonzero
$K$-multiple of $t$.

A cycle-basis change acts on both rows by the same element of
$\GL_2(\Qp)$. Rationality over $\Qp$ and algebraicity over $\Qp$
of a projective row are therefore independent of the chosen
$\Qp$-basis of $V$, with the de Rham class fixed.
Replacing $\xi$ by $\xi+c\omega$ changes the second projective point
and changes the de Rham class whose period row enters the
rationality condition.

For $\eta\in\HdR(E/K)$, write
\[
 \Ann_V(\eta)=\{v\in V:\langle v,\eta\rangle=0\}.
\]

\begin{lemma}[Rational rows and stable lines]\label{lem:rows}
For nonzero $\eta\in\HdR(E/K)$, its period row is in $\Pj^1(\Qp)$
if and only if $\Ann_V(\eta)$ has dimension one over $\Qp$.
That annihilator is $G_K$-stable. Conversely, every stable line
$M\subset V$ is the annihilator of a unique $K$-line
$U_M\subset\HdR(E/K)$, whose projective period row is rational.
\end{lemma}
\begin{proof}
A nonzero row $(A,B)$ has a nonzero $\Qp$-linear relation precisely
when $[A:B]\in\Pj^1(\Qp)$. Its kernel cannot have dimension two
because comparison is perfect. Equivariance makes the kernel stable.
For the converse, $\HdR(E/K)=\DdR(V^*)$.
The annihilator of $M$ in $V^*$ is $(V/M)^*$, a one-dimensional
de Rham representation. Exactness of $\DdR$ gives a $K$-line
in $\HdR(E/K)$; comparison identifies it with the annihilator
of $M\otimes\BdR$. This proves existence and uniqueness.
\end{proof}

\begin{lemma}[Two row conditions force splitting]\label{lem:splitting}
If $[H_1:H_2]\in\Pj^1(\Qp)$ and
$[W_1:W_2]\in\Pj^1(\Qpb)$, then $V$ becomes a direct sum
of two one-dimensional $\Qp$-representations over a finite extension.
\end{lemma}
\begin{proof}
Choose $e_1$ on the stable annihilator of $\xi$.
Then $H_1=0$, $H_2\ne0$, and $W_1\ne0$. Put
$\alpha=W_2/W_1\in\Qpb$, and choose finite $L/K$ fixing $\alpha$.
For $\sigma\in G_L$ write
\[
 \sigma(e_1)=a_\sigma e_1,\qquad
 \sigma(e_2)=b_\sigma e_1+d_\sigma e_2,
\]
with $a_\sigma,d_\sigma\in\Qp^\times$ and $b_\sigma\in\Qp$.
Since $\omega$ is defined over $K$, equivariance gives
\[
 \sigma(W_1)=a_\sigma W_1,\qquad
 \sigma(W_2)=b_\sigma W_1+d_\sigma W_2.
\]
For $\sigma\in G_L$, the equality $W_2=\alpha W_1$ and
$\sigma(\alpha)=\alpha$ therefore imply
$a_\sigma\alpha=b_\sigma+d_\sigma\alpha$. Hence
\begin{equation}
 b_\sigma=(a_\sigma-d_\sigma)\alpha. \label{eq:splitting-coeff}
\end{equation}
If some $a_\sigma\ne d_\sigma$, this forces $\alpha\in\Qp$,
and $\Qp(e_2-\alpha e_1)$ is a stable complement to $\Qp e_1$.
Otherwise every $b_\sigma=0$, and $\Qp e_2$ is a stable complement.
\end{proof}
For $q\in K^\times$ with $0<|q|_p<1$, let $E_q/K$ denote
the Tate elliptic curve, with rigid analytic uniformization
\[
 E_q^{\mathrm{an}}
 \simeq (\mathbb G_{m,K})^{\mathrm{an}}/q^{\Z};
\]
see \cite[Appendix~C, Theorem~14.1]{SilvermanAEC}.
The same parameter defines the Kummer $1$-motive
\[
 M_q=[\,\Z\xrightarrow{\,1\mapsto q\,}\Gm\,].
\]
Tate uniformization identifies its $p$-adic Tate realization
$T_p(M_q)$ with $T_p(E_q)$, and hence gives an exact sequence
\[
 0\longrightarrow\Z_p(1)\longrightarrow T_p(E_q)
 \longrightarrow\Z_p\longrightarrow0.
\]
Its extension class is $\kappa_K(q)$, where
\[
 \kappa_K:K^\times\longrightarrow H^1(G_K,\Z_p(1))
\]
is the Kummer map induced by the $p^n$-power maps on $\Gm$;
see \cite[\S II.4]{Berger}.
\begin{lemma}[Nonsplitting of a Tate extension]\label{lem:kummer}
For $q\in K^\times$ with $0<|q|_p<1$, the sequence
\[
 0\longrightarrow\Qp(1)\longrightarrow V_p(E_q)
 \longrightarrow\Qp\longrightarrow0
\]
is nonsplit after every finite extension of $K$.
\end{lemma}
\begin{proof}
The rational Tate-module extension is classified by the image
of $\kappa_K(q)$ under the natural map
\[
 H^1(G_K,\Z_p(1))
 \longrightarrow H^1(G_K,\Qp(1)).
\]
For finite $L/K$, its integral class is the image of $q$ in
$\varprojlim_nL^\times/(L^\times)^{p^n}$.
Valuation maps this completion to $\Z_p$ and sends $q$ to
$v_L(q)>0$. This remains nonzero after tensoring with $\Qp$,
so the rational extension class does not vanish.
The representation and its Kummer description are recalled
in~\cite[\S II.4]{Berger}.
\end{proof}

\begin{theorem}[Local quasi-period implication]\label{thm:local-A}
Suppose $\End_{\overline K}(E)=\Z$.
For every basis $\omega,\xi$ of $\HdR(E/K)$ with $\omega$ holomorphic,
\begin{equation}
 [H_1:H_2]\in\Pj^1(\Qp)\quad\Longrightarrow\quad
 [W_1:W_2]\notin\Pj^1(\Qpb).
 \tag{$A_p$}\label{eq:Ap}
\end{equation}
Thus a rational quasi-period ratio forces the period ratio to be
transcendental over $\Qp$. This is the local form of Question A.
The analogue of Kumar's zero reformulation appears on the ordinary
locus in Theorem~\ref{thm:ordinary-condition}.
\end{theorem}
\begin{proof}
Suppose that both row conditions hold. By
Lemma~\ref{lem:splitting}, $V_p(E)$ is a direct sum of two
one-dimensional representations after a finite extension.
The semistable reduction theorem, followed if necessary by an
extension splitting multiplicative reduction, gives a finite
extension $L/K$ over which $E_L$ has good or split
multiplicative reduction; see
\cite[Chapter~VII, \S5, and Appendix~C, \S14]{SilvermanAEC}.
The decomposition of the representation persists upon
restriction to $G_L$.

In the good-reduction case, Serre's corollary
\cite[Chapter~IV, Appendix~A.2.3]{Serre} applies: an abelian
variety over a local field with finite residue field whose
rational Tate representation is a sum of $\Qp$-characters
is isogenous to a product of CM abelian varieties.
In dimension one this forces geometric CM, contradicting
$\End_{\overline K}(E)=\Z$.

In the split multiplicative case, Tate's uniformization theorem
\cite[Appendix~C, Theorem~14.1]{SilvermanAEC} gives
$E_L\simeq E_q$ over $L$ for some $q\in L^\times$ with
$0<|q|_p<1$. Thus
$V_p(E)|_{G_L}\simeq V_p(E_q)$. A direct sum of characters
is semisimple, so its cyclotomic subrepresentation has a
$G_L$-stable complement. This would split the Tate extension,
contrary to Lemma~\ref{lem:kummer}. Both reduction cases
therefore exclude simultaneous satisfaction of the two row
conditions.
\end{proof}
Serre's finite-residue-field hypothesis applies because $K$ is
a finite extension of $\Qp$. The proof uses potential splitting
of the rational Tate representation; it makes no assumption that
every local Tate-module endomorphism comes from an elliptic
endomorphism.

\begin{corollary}\label{cor:irreducible}
If $V_p(E)$ is irreducible, no nonzero de Rham class has a rational
projective period row. This holds for good supersingular reduction,
even after any finite extension.
\end{corollary}
\begin{proof}
The first assertion follows from Lemma~\ref{lem:rows}.
For good supersingular reduction, the crystalline realization
of $V^*$ is isoclinic of slope $1/2$. A line in $V$ would dualize
to a one-dimensional crystalline quotient of $V^*$, hence to
a rank-one quotient isocrystal of slope $1/2$. A rank-one
isocrystal over the maximal unramified subfield $K_0$ of $K$
has integral slope. Indeed, write Frobenius as $a\sigma$,
where $a\in K_0^\times$ and $\sigma$ is the arithmetic
Frobenius of the unramified extension $K_0/\Qp$.
If $f=[K_0:\Qp]$, then its $f$-fold iterate is
$N_{K_0/\Qp}(a)$, so the slope is
\[
 \frac1f v_p\bigl(N_{K_0/\Qp}(a)\bigr)
 =v_p(a)\in\Z.
\]
It therefore cannot equal $1/2$. The same calculation holds after any finite
extension, where supersingularity is preserved.
\end{proof}

\subsection{Period relations modulo powers of \texorpdfstring{$t$}{t}}
\label{subsec:finite-truncation}
The Hodge filtration bounds the power of $t$ needed to detect
algebraic linear period relations. For each integer $n\ge1$, put
\[
 B_n=\BdRp/t^n\BdRp,\qquad
 q_n:\BdRp\longrightarrow B_n,\qquad
 \Fil^n\BdR=t^n\BdRp.
\]
Thus $q_n(x)$ records a period $x\in\BdRp$ up to addition
of an element of $\Fil^n\BdR$. This quotient operation is
the $t$-adic truncation used here. In particular, the comparison
pairing modulo $t^2$ is the composite of
\eqref{eq:pairing} with $q_2$; its value at $(v,\eta)$ is
$q_2(\langle v,\eta\rangle)\in B_2$.

The following lemma applies to any de Rham representation with
a bounded nonnegative filtration. We use strict exactness and
tensor compatibility of the filtered de Rham functor
\cite[Proposition~6.3.3]{BC}; finite base change is compatible
with this filtration by \cite[Proposition~6.3.8]{BC}.

\begin{lemma}[A finite filtration bound]\label{lem:finite-filtration}
Let $R$ be a de Rham representation of $G_K$, and let $h\in\Z_{\ge0}$.
Suppose
\[
 \Fil^0\DdR(R^*)=\DdR(R^*),\qquad
 \Fil^{h+1}\DdR(R^*)=0.
\]
For $v\in R$ and $\eta\in\DdR(R^*)$, one has
\[
 \langle v,\eta\rangle\in t^{h+1}\BdRp
 \quad\Longleftrightarrow\quad \langle v,\eta\rangle=0.
\]
\end{lemma}
\begin{proof}
Let $W=\operatorname{span}_{\Qp}\{\sigma v:\sigma\in G_K\}$.
This is a de Rham subrepresentation of $R$. If the displayed
period belongs to $t^{h+1}\BdRp$, equivariance gives
\[
 \langle\sigma v,\eta\rangle
   =\sigma\langle v,\eta\rangle\in t^{h+1}\BdRp
 \qquad(\sigma\in G_K).
\]
Thus the restriction $\eta_W\in\DdR(W^*)$ takes values in
$t^{h+1}\BdRp$. By the induced filtration on
$\BdR\otimes W^*=\Hom_{\Qp}(W,\BdR)$, this says
$\eta_W\in\Fil^{h+1}\DdR(W^*)$. Strict exactness identifies
the filtration on $\DdR(W^*)$ with the quotient filtration
from $\DdR(R^*)$. Its $(h+1)$st step is zero, so $\eta_W=0$.
Evaluation at $v$ proves the forward implication. The reverse
implication is immediate.
\end{proof}

\begin{proposition}[Linear relations detected in $B_2$]
\label{prop:linear-truncation}
For the full elliptic comparison matrix \eqref{eq:matrix},
the map $q_2$ is injective on
\[
 \mathcal P_{\le1}(E)
 =\operatorname{span}_{\Qpb}\{1,W_1,W_2,H_1,H_2\}
 \subset\BdRp.
\]
Consequently every linear relation over $\Qpb$ among these
comparison periods in \eqref{eq:matrix} holds in the full period ring if
and only if it holds after reduction to $B_2$.
\end{proposition}
\begin{proof}
The dual de Rham realization $D=\HdR(E/K)=\DdR(V^*)$
satisfies $\Fil^0D=D$ and $\Fil^2D=0$.
Given a finite linear combination of the displayed generators,
choose finite $L/K$ containing its coefficients.
Absorb each coefficient into the corresponding de Rham class
over $L$. The sum is a single comparison pairing on a finite
direct sum of copies of $V$; a copy of the trivial
representation accounts for the constant term. The dual de Rham
realization of this sum still has $\Fil^2=0$ and $\Fil^0$
equal to the whole space. Lemma~\ref{lem:finite-filtration}
with $h=1$ proves the claim.
\end{proof}

For good-reduction $1$-motives, the comparison of vanishing in
the crystalline pairing and the pairing modulo $t^2$ is recorded in
\cite[Lemma~3.1 and Corollary~3.2]{MS}.
Proposition~\ref{prop:linear-truncation} formulates the result
for all algebraic linear combinations of elliptic comparison
entries. Its filtered proof applies to every reduction type.
The same argument gives the following extension.

\begin{corollary}[Relations of bounded degree]
\label{cor:bounded-degree-truncation}
Let $d\in\Z_{\ge0}$, and let $\mathcal P_{\le d}(E)$ be the
$\Qpb$-span of all monomials of total degree at most $d$ in
the four entries of $P$. Reduction $q_{d+1}$ is injective
on $\mathcal P_{\le d}(E)$.
\end{corollary}
\begin{proof}
A degree-$r$ monomial is a comparison pairing for
$V^{\otimes r}$, with degree zero corresponding to the trivial
representation. Its dual de Rham filtration is concentrated in
degrees $0,\ldots,r$. After putting the coefficients in a finite
extension, a finite sum is a pairing on a direct sum of these
tensor powers. Apply Lemma~\ref{lem:finite-filtration} with
$h=d$.
\end{proof}

The linear argument applies to $H^1$ of any abelian variety
over $K$, and more generally to any de Rham representation
whose dual filtration is concentrated in degrees zero and one.
The filtration range determines the order of truncation.
The bound in Corollary~\ref{cor:bounded-degree-truncation} is
sharp for every elliptic curve and every $d\ge1$. The
holomorphic row has a nonzero entry $W_i$ by perfection of the
comparison. The polynomial $X_{1i}^d$ takes the nonzero value
$W_i^d\in t^d\BdRp$, whose image vanishes in $B_d$ and hence
in every $B_n$ with $n\le d$.

The bound measures the $t$-adic filtration. In $B_1$ the entire
holomorphic row vanishes. In $B_2$ the determinant is a nonzero
multiple of the nilpotent class of $t$. The truncation results
test polynomial equations in the entries of $q_2(P)$ as a
matrix over $B_2$.

\subsection{Holomorphic relations in the degree-one comparison}
\label{subsec:holomorphic-graded}
Define
\[
 a_i=\theta(W_i/t)\in\Cp\qquad(i=1,2).
\]
These are the degree-one holomorphic comparison periods after
trivializing the Tate twist by the class of $t$.

\begin{proposition}[Full and degree-one holomorphic relations]
\label{prop:holomorphic-graded}
The row $(a_1,a_2)$ is nonzero. For $c_1,c_2\in\Qpb$,
\[
 c_1W_1+c_2W_2=0
 \quad\Longleftrightarrow\quad c_1a_1+c_2a_2=0.
\]
In particular, $[W_1:W_2]\in\Pj^1(\Qpb)$ if and only if
$[a_1:a_2]\in\Pj^1(\Qpb)$. Whenever these conditions hold,
the two points are the images of the same point of
$\Pj^1(\Qpb)$ in their respective comparison fields.
The same assertions hold with $\Qp$ in place of $\Qpb$.
\end{proposition}
\begin{proof}
The right-hand relation means
$c_1W_1+c_2W_2\in t^2\BdRp$, because $\theta$ restricts
to the fixed embedding of $\Qpb$ in $\Cp$.
Proposition~\ref{prop:linear-truncation} makes this equivalent
to full vanishing. If both $a_i$ were zero, that proposition
would give $W_1=W_2=0$, contradicting comparison. The
projective assertions follow by taking the one-dimensional
space of relations of either nonzero row.
\end{proof}

The intrinsic map behind this row is
\begin{equation}
 \Phi_E:V_p(E)\longrightarrow
   \Hom_K(\Fil^1\HdR(E/K),\Cp(1))
   =\Lie(E)\otimes_K\Cp(1),\qquad
 \Phi_E(v)(\omega)=\operatorname{gr}^1\langle v,\omega\rangle.
 \label{eq:Fontaine-map}
\end{equation}
Here $\Cp(1)=t\BdRp/t^2\BdRp$ as a $G_K$-module.
After choosing $\omega$ and the class of $t$, the coordinate
row of $\Phi_E$ is $(a_1,a_2)$.
For good reduction this is Fontaine's map used in
\cite[Section~3.2 and equation~(3.19)]{MS}.
This fixes the target line, Tate twist, and normalization.
Its relation with the canonical Hodge--Tate filtration requires
distinguishing an inclusion from a projection. The filtration
sequence for $E[p^\infty]$, using the principal polarization, is
\begin{equation}
 0\longrightarrow\Lie(E)\otimes_K\Cp(1)
 \xrightarrow{\ \iota_{\mathrm{HT}}\ }V\otimes_{\Qp}\Cp
 \longrightarrow\Lie(E^\vee)^*\otimes_K\Cp\longrightarrow0.
 \label{eq:HT-filtration}
\end{equation}
The image $L_1$ of $\iota_{\mathrm{HT}}$ is the filtration point
studied by Howe~\cite[\S1.2 and Theorem~1.2]{Howe}. For the model
over $K$, the $G_K$-equivariant Hodge--Tate decomposition supplies
a complementary weight-zero line $L_0$, and
\[
 V\otimes_{\Qp}\Cp=L_1\oplus L_0.
\]
The scalar extension of $\Phi_E$ projects onto $L_1$ and then
identifies $L_1$ with $\Lie(E)\otimes_K\Cp(1)$. Its covector
$(a_1,a_2)$ describes $L_0$ by its kernel. The filtration is
defined for the $p$-divisible group over $\mathcal O_{\Cp}$;
the equivariant complement comes from the model over $K$.

For ordinary reduction, $L_1=V^+\otimes_{\Qp}\Cp$, where $V^+$
is the connected Tate line, so the included line is rational
over $\Qp$. Howe--Klevdal
\cite[Remark~1.1.4]{HoweKlevdal} describe the corresponding
gap between filtration and grading. They relate algebraicity
of the grading to torsion of the Serre--Tate parameter, using
Serre's orbit calculation for the Hodge--Tate cocharacter
\cite[\S1.5, Theorem~3]{SerreHT}. In these ordinary coordinates
algebraicity of the grading amounts to algebraicity of $L_0$.
Corollary~\ref{cor:holomorphic-CM} derives the criterion for the
\emph{full} period row from local monodromy, and
Proposition~\ref{prop:holomorphic-graded} identifies its
projective algebraicity with that of the projection row.

The bases for the corresponding CM statements matter as well.
Howe's CM means that
$\End_{\mathcal O_{\Cp}}(G)\otimes_{\Z_p}\Qp$ contains a
commutative semisimple algebra of dimension $\operatorname{ht}G$.
In particular, ordinary extensions split over $\mathcal O_{\Cp}$
and have CM there; that splitting need not descend to any finite
extension of $K$ \cite[Remark~1.3]{Howe}. Geometric CM here means
an extra endomorphism of $E_{\overline K}$. The formal CM used in
Table~\ref{tab:local} concerns endomorphisms of the height-two formal
group after a finite extension of the local base. Our criterion
The monodromy argument below makes these distinctions explicit
in the relation ideals of the two comparison matrices.

\section{Local monodromy and the full period field}\label{sec:monodromy}
\subsection{The full comparison matrix}
The Galois action on a full $\BdR$ comparison matrix gives a lower
bound for the transcendence degree of its entries over $\Qpb$.
The bound is the dimension of the Zariski closure of the Galois image.

\begin{lemma}[Galois-orbit lower bound]\label{lem:orbit}
Let $P\in\GL_n(\BdR)$ satisfy $\sigma(P)=P\rho(\sigma)$ for
a continuous representation $\rho:G_K\to\GL_n(\Qp)$.
If $G$ is the Zariski closure of $\rho(G_K)$ over $\Qp$, then
\[
 \trdeg_{\Qpb}\Qpb(P)\geq\dim G.
\]
\end{lemma}
\begin{proof}
Let $I$ be the ideal of polynomial relations on $P$ over $\Qpb$,
using the coordinate ring of $\GL_n$. For any $f\in I$, choose
finite $L/K$ containing its coefficients. For $\sigma\in G_L$,
\[
 f(P\rho(\sigma))=\sigma(f(P))=0.
\]
The closure of $\rho(G_L)$ contains the identity component $G^\circ$.
Zariski density still holds after extension of scalars to $\BdR$:
express the finitely many coefficients of a polynomial over $\BdR$
in a $\Qp$-linearly independent finite set to reduce to polynomials
over $\Qp$. It follows that $f(Pg)=0$ for every
$g\in G^\circ_{\BdR}$. Thus the closed subvariety of $\GL_n$ defined by $I$, after
base change to $\BdR$, contains $P G^\circ_{\BdR}$. Right multiplication on an invertible
matrix is free, so this orbit has dimension $\dim G$.
Dimension is preserved by field extension, giving the lower bound.
\end{proof}

By formal CM in the supersingular rows of Table~\ref{tab:local}
we mean CM of the height-two formal group over $\mathcal O_{K'}$
for some finite extension $K'/K$, in Serre's local sense.
Equivalently in this classification, its Tate representation over
$G_{K'}$ has a quadratic field over $\Qp$ of commuting endomorphisms.
Geometric CM concerns $\End_{\overline K}(E)$; Howe's CM over
$\mathcal O_{\Cp}$ was specified after~\eqref{eq:HT-filtration}.

\begin{theorem}[Transcendence degree of the full local matrix]
\label{thm:dimensions}
For an elliptic curve over finite $K/\Qp$ and a period matrix
\eqref{eq:matrix}, the value
$d_p(E)=\trdeg_{\Qpb}\Qpb(P)$ is given by
Table~\ref{tab:local}. Here $\Qpb(P)$ denotes the field generated
by the four entries. Reduction types are taken after a finite
extension, and $G^\circ$ is the connected algebraic Galois group.
\end{theorem}
\begin{table}[htbp]
\centering
\(
\begin{array}{@{}lcc@{}}
\toprule
\text{Reduction and endomorphism type}&G^\circ&d_p(E)\\
\midrule
\text{Split multiplicative}&\Gm\ltimes\Ga&2\\
\text{Good ordinary, no geometric CM}&\text{Borel in }\GL_2&3\\
\text{Good ordinary, geometric CM}&\text{split torus}&2\\
\text{Good supersingular, no formal CM}&\GL_2&4\\
\text{Good supersingular, formal CM}&\text{nonsplit torus}&2\\
\bottomrule
\end{array}
\)
\caption{Full local de Rham period fields and connected monodromy.}
\label{tab:local}
\end{table}
\begin{proof}
Carry the chosen bases through de Rham base change when
passing to a finite extension of $K$. The comparison entries
then stay fixed in $\BdR$, and the coefficient field remains
$\Qpb$. Invertible algebraic changes of de Rham or cycle
coordinates preserve their generated field and transcendence degree.

Serre's local description gives the connected groups in the
middle column: see \cite[Chapter~IV, Appendix~A.1]{Serre}
for multiplicative reduction, Appendix~A.2.2 for supersingular
reduction, and Appendix~A.2.4 for ordinary reduction.
These descriptions determine the connected algebraic groups
from their local Galois images and Lie algebras.
In the multiplicative case the group has matrices
$\left(\begin{smallmatrix}a&b\\0&1\end{smallmatrix}\right)$.
In the ordinary non-CM case it is the full triangular group.
In the two CM cases it is a torus of dimension two, split in the
ordinary case and nonsplit in the supersingular case.
Without formal CM, the supersingular image has full Lie algebra
$\mathfrak{gl}_2$. Lemma~\ref{lem:orbit} supplies the displayed
lower bounds.

For the upper bounds, four entries give $d_p(E)\leq4$.
A stable $\Qp$-line in the ordinary case, together with
Lemma~\ref{lem:rows}, gives algebraic coordinate choices with
\[
 P=\begin{pmatrix}A&C\\0&B\end{pmatrix},
\]
so $d_p(E)\leq3$.
In the multiplicative case the quotient character is trivial.
The second row can therefore be normalized to $(0,1)$:
its nonzero entry is $G_K$-fixed and hence belongs to $K$.
This gives $d_p(E)\leq2$.

In either torus case, after finite extension there is a nonscalar
commuting endomorphism $J$ of $V$ with separable quadratic minimal
polynomial. Applying $\DdR$ to this Galois-equivariant map
gives a filtered endomorphism $J_{\mathrm{dR}}$ over that
finite extension of $K$, and comparison intertwines the two.
This construction uses the endomorphism of the representation
itself, including in the formal CM case.
Extend the coefficient field to split the minimal polynomial.
Since $J_{\mathrm{dR}}$ preserves the Hodge line, keep $\omega$
as a first eigenvector and take the other eigenvector as the
second de Rham basis element. Order the two cycle eigenvectors
to match. Comparison then makes $P$ diagonal, giving $d_p(E)\leq2$.
All lower and upper bounds agree.
\end{proof}

This theorem concerns full comparison periods over $\Qpb$.
It is not a theorem about algebraic independence over $\Q$ of
$\Cp$-valued Coleman integrals or of crystalline period invariants.
Its proof uses the Galois action on $\BdR$, which is absent from
the corresponding complex comparison field.

\begin{corollary}[The ideal of local algebraic period relations]
\label{cor:local-relation-ideal}
Let $S=\Qpb[X_{11},X_{12},X_{21},X_{22}]$, and let
$I_E$ be the kernel of evaluation at the full matrix $P$.
Then $I_E$ is generated by polynomials of degree at most one.
More precisely, after invertible linear changes of de Rham
and cycle coordinates over $\Qpb$, the relation ideals and
the numbers of algebraically independent entries are
\[
\begin{array}{@{}lcc@{}}
 \text{Monodromy type}&I_E&d_p(E)\\\midrule
 \GL_2&0&4\\
 \text{ordinary non-CM Borel}&(X_{21})&3\\
 \text{either CM torus}&(X_{12},X_{21})&2\\
 \text{multiplicative}&(X_{21},X_{22}-1)&2.
\end{array}
\]
If $S_{\le1}$ is the space of affine linear polynomials and
$\overline{\mathrm{ev}}_{P,2}:S\to B_2$ is evaluation at
$q_2(P)$, then
\[
 I_E=\bigl(\ker(\overline{\mathrm{ev}}_{P,2})
                   \cap S_{\le1}\bigr),
\]
where parentheses denote the ideal generated in $S$.
\end{corollary}
\begin{proof}
The coordinate choices in the proof of
Theorem~\ref{thm:dimensions} put $P$ respectively in full,
upper triangular, diagonal, or upper triangular form with
bottom row $(0,1)$. In each case the number of entries left
free is exactly $d_p(E)$ in Table~\ref{tab:local}.
Those entries are therefore algebraically independent over
$\Qpb$. The displayed equations generate the full kernel of
evaluation. Invertible linear coordinate changes preserve
generation by affine linear polynomials.
Finally, Proposition~\ref{prop:linear-truncation} identifies
the affine linear kernels of evaluation at $P$ and at $q_2(P)$.
\end{proof}

This describes algebraic relations over the local algebraic
coefficient field, with Serre's monodromy classification as
input. Higher-degree equations can vanish in $B_2$ through
nilpotence: for example, $\det(q_2(P))^2=0$ whereas
$\det(P)^2\ne0$. The last assertion recovers the full relation
ideal from its affine linear generators; testing a polynomial
directly is governed by
Corollary~\ref{cor:bounded-degree-truncation}.

\subsection{Relations of the normalized comparison matrix}
\label{subsec:graded-matrix}
Put $\delta=\det P/t\in K^\times$. Divide the holomorphic row
of $P$ by $t$ and reduce its entries and the second row using
$\theta:\BdRp\to\Cp$. The resulting matrix is
\begin{equation}
 P_{\mathrm{gr}}=
 \begin{pmatrix}
   \theta(W_1/t)&\theta(W_2/t)\\
   \theta(H_1)&\theta(H_2)
 \end{pmatrix},
 \qquad \det P_{\mathrm{gr}}=\delta .
 \label{eq:graded-matrix}
\end{equation}
This matrix records the degree-one and degree-zero comparison
rows. If $\omega$ is multiplied by $a\ne0$ and $\xi$ is
replaced by $b\xi+c\omega$ with $b\ne0$, the two rows of
$P_{\mathrm{gr}}$ are multiplied by $a$ and $b$, respectively.
The term $c\omega$ has zero reduction in degree zero.
Algebraic changes of cycle coordinates act on the right.
The normal forms used below therefore use the same choices
as the full matrix, keeping the first de Rham row holomorphic.

\begin{theorem}[Relation ideal of the normalized comparison matrix]
\label{thm:graded-relation-ideal}
Let $I_{\mathrm{gr}}$ be the kernel of evaluation at
$P_{\mathrm{gr}}$ in $\Qpb[X_{11},X_{12},X_{21},X_{22}]$.
In the normal forms of Corollary~\ref{cor:local-relation-ideal},
the ideal and transcendence degree are as follows.
\begin{equation}
\label{eq:graded-ideals}
\begin{array}{@{}lcc@{}}
\toprule
\text{Monodromy type} & I_{\mathrm{gr}} &
       \trdeg_{\Qpb}\Qpb(P_{\mathrm{gr}})\\
\midrule
\GL_2 & (\det X-\delta)&3\\
\text{ordinary non-CM Borel}&(X_{21},\det X-\delta)&2\\
\text{split or nonsplit torus}&(X_{12},X_{21},\det X-\delta)&1\\
\text{multiplicative}&(X_{21},X_{22}-1,X_{11}-\delta)&1\\
\bottomrule
\end{array}
\end{equation}
In particular, the transcendence degree is $d_p(E)-1$.
\end{theorem}
\begin{proof}
Enlarge $K$ to fix any algebraic basis-change coefficients.
Let $F=\Qpb$ and write $g=\rho(\sigma)$ for the matrix of
$\sigma\in G_K$ on the cycle space. The Weil pairing gives
$\det g=\chi(\sigma)$, where $\chi$ is the cyclotomic character.
The full matrix satisfies $\sigma(P)=Pg$ and
$\sigma(t)=\chi(\sigma)t$, so equivariance of $\theta$ gives
\begin{equation}
 \sigma(P_{\mathrm{gr}})=
 \operatorname{diag}\bigl((\det g)^{-1},1\bigr)P_{\mathrm{gr}}g .
 \label{eq:graded-action}
\end{equation}
The same formula defines an algebraic right action of
$G^\circ_{\Cp}$ on invertible matrices over $\Cp$.

For each polynomial relation over $F$, choose a finite extension
fixing its coefficients. Evaluating its Galois conjugates and
using the Zariski-density argument of Lemma~\ref{lem:orbit}
shows that the base change of its zero set to $\Cp$ contains
the orbit of $P_{\mathrm{gr}}$ under $G^\circ_{\Cp}$.
This holds for every relation in $I_{\mathrm{gr}}$. A group
element stabilizing $P_{\mathrm{gr}}$ must satisfy
\[
 g=P_{\mathrm{gr}}^{-1}\operatorname{diag}(\det g,1)P_{\mathrm{gr}},
\]
so the stabilizer is contained in a one-dimensional torus.
The orbit dimension therefore gives
\[
 \trdeg_F F(P_{\mathrm{gr}})\ge\dim G^\circ-1=d_p(E)-1.
\]

The listed relations hold by~\eqref{eq:graded-matrix} and the
full normal forms. Each of the displayed ideals is prime. The
fiber $\det X=\delta$ is isomorphic to $\operatorname{SL}_2$.
With $X_{21}=0$ the same determinant equation gives
$X_{11}X_{22}=\delta$ and leaves $X_{12}$ free. With both
off-diagonal entries zero it gives $X_{11}X_{22}=\delta$.
In the last case $X_{12}$ is free and the remaining entries are fixed.
The four quotients have dimensions $3,2,1,1$, respectively.
In each case the quotient attains the lower bound, so the displayed
prime ideal is the entire relation ideal.
\end{proof}

Serre's orbit formula
\cite[\S1.5, Theorem~3]{SerreHT} concerns the field of
the Hodge--Tate grading, obtained from conjugation of its
cocharacter. Equation~\eqref{eq:graded-action} retains the
normalized entries of both rows. For example, in the ordinary
case the full graded matrix keeps a unit-root period as well
as the coordinate of the projection line.
Corollary~\ref{cor:ordinary-zero-periods} gives the resulting
joint independence.

If $G^\circ=\GL_2$, both projective row coordinates of
$P_{\mathrm{gr}}$ are algebraically independent over $\Qpb$.
Indeed, the determinant fiber in~\eqref{eq:graded-ideals}
maps dominantly to pairs of distinct points of
$\Pj^1\times\Pj^1$: scale two representatives of any such pair
so that their determinant equals $\delta$. The same
independence holds for the full row coordinates by
Corollary~\ref{cor:local-relation-ideal}.

\subsection{Holomorphic period points and the CM obstruction}

\begin{corollary}[Algebraic holomorphic period points]
\label{cor:holomorphic-CM}
For an elliptic curve over finite $K/\Qp$, the full holomorphic
point $[W_1:W_2]$ is algebraic over $\Qp$ precisely in the
following cases:
\begin{enumerate}
\item potentially good ordinary reduction with geometric CM;
\item potentially good supersingular reduction with formal CM.
\end{enumerate}
In the ordinary case the point belongs to $\Pj^1(\Qp)$.
In the supersingular case it belongs to
$\Pj^1(F)\setminus\Pj^1(\Qp)$, where $F/\Qp$ is the
quadratic local CM field.
The degree-one point $[a_1:a_2]$ has the same algebraicity
criterion and the same field of definition.
\end{corollary}
\begin{proof}
If the full point is algebraic, choose finite $L/K$ fixing
its projective coordinates. Equivariance makes its row line
invariant under the Galois image over $L$, hence under its
connected Zariski closure over $\Qpb$.
The group $\GL_2$ has no invariant dual line, excluding the
supersingular case without formal CM in Table~\ref{tab:local}.

In the ordinary case, let $V^+$ be the Tate line of the
connected height-one part of the ordinary $p$-divisible group.
In the split multiplicative case, let $V^+=\Qp(1)$ be the
cyclotomic subrepresentation of the Tate extension. The
connected groups in the two non-CM cases have a unique
invariant dual line, namely the annihilator of $V^+$.

The dual of $V^+$ has its de Rham filtration concentrated in
degree one. Strict exactness applied to the restriction map
therefore gives a surjection
\[
 K\omega=\Fil^1\DdR(V^*)\longrightarrow
 \Fil^1\DdR((V^+)^*)=\DdR((V^+)^*).
\]
Both spaces have dimension one, so $\omega$ restricts
nontrivially to $V^+$. Its period row cannot be the invariant
dual line. This excludes both cases.

In either torus case, choose a nonscalar commuting endomorphism
$J$ with separable quadratic minimal polynomial, as in the
proof of Theorem~\ref{thm:dimensions}. Its filtered de Rham
endomorphism preserves the one-dimensional Hodge line.
Comparison makes the holomorphic row an eigenvector of the
$\Qp$-matrix of $J$. In the split case both eigenlines are
defined over $\Qp$. In the nonsplit case they are defined over
its quadratic splitting field $F$, and neither is defined over
$\Qp$. This proves the converse and the stated fields.
Proposition~\ref{prop:holomorphic-graded} gives the degree-one
assertions.
\end{proof}
This criterion concerns the projection covector of
Subsection~\ref{subsec:holomorphic-graded}, independently of the
second-kind class. It distinguishes formal CM over a finite local
extension from geometric CM. Howe's theorem concerns the included
line in~\eqref{eq:HT-filtration}; the ordinary example there shows
why it cannot be substituted for the present monodromy argument.

\begin{samepage}
\begin{corollary}[The exact local obstruction]\label{cor:exact-obstruction}
Let $\omega,\xi$ be as in Theorem~\ref{thm:local-A}, without the
non-CM hypothesis. Both
\[
 [H_1:H_2]\in\Pj^1(\Qp),\qquad
 [W_1:W_2]\in\Pj^1(\Qpb)
\]
hold if and only if there is a finite extension $K'/K$ over which
$E$ has good ordinary reduction with geometric CM and $K'\xi$
is the complementary CM de Rham eigenspace, equivalently the
unit-root line.
In that case the holomorphic period row is already $\Qp$-rational.
\end{corollary}
\end{samepage}
\begin{proof}
The two row conditions imply potential splitting by
Lemma~\ref{lem:splitting}. Multiplicative reduction is excluded by
Lemma~\ref{lem:kummer}, and supersingular reduction by
Corollary~\ref{cor:irreducible}. Serre's corollary forces CM in the
remaining ordinary case. Over a field of definition of the CM action
the representation splits into two distinct characters. Its only
stable lines are the two eigenspaces. Their annihilator lines in
cohomology are the Hodge line $K'\omega$ and its complementary CM
line, which is the unit-root line. Since $\xi$ is independent of
$\omega$, the rational row for $\xi$ must be the second one.
Conversely, the two CM eigenspaces and their annihilators give
rational rows by Lemma~\ref{lem:rows}.
These projective-row conditions remain unchanged after extending $K$, with the comparison matrix and embeddings fixed.
\end{proof}

\section{The Tate curve: formulas and independence}\label{sec:tate}
Let $q\in K^\times$ satisfy $0<|q|_p<1$. Retain the
cyclotomic system $\varepsilon$ fixed in Section~\ref{sec:local},
choose compatible roots
$\widetilde q=(q,q^{1/p},q^{1/p^2},\ldots)$, and put
\begin{equation}
 t=\log[\varepsilon],\qquad
 \ell_q=\log\bigl([\widetilde q]/q\bigr)\in t\BdRp,\qquad
 \tau_q=\ell_q/t.
 \label{eq:tate-periods}
\end{equation}
The quotient $[\widetilde q]/q$ belongs to $1+\ker\theta$.
The logarithm here is the convergent $\ker\theta$-adic logarithm.
In particular $\ell_q$ is not the scalar $\log_p(q)$.
In Berger's notation it is $u-\log_p(q)$; the filtration in
\cite[\S II.4.3]{Berger} gives precisely this distinction.

The Tate basis $e,f$ comes from $\varepsilon$ and $\widetilde q$.
Its action is
\begin{equation}
 \sigma(e)=\chi(\sigma)e,\quad
 \sigma(f)=f+c_q(\sigma)e,\quad
 \sigma(\ell_q)=\ell_q+c_q(\sigma)t,\quad
 \sigma(t)=\chi(\sigma)t .
 \label{eq:tate-action}
\end{equation}
The Kummer cocycle $c_q(\sigma)\in\Z_p$ is defined by the
compatible roots-of-unity quotients
$\sigma(q^{1/p^n})/q^{1/p^n}=\zeta_{p^n}^{c_q(\sigma)}$.
In particular
$c_q(\sigma\tau)=c_q(\sigma)+\chi(\sigma)c_q(\tau)$.

\subsection{The standard short model}
Write $s_1(q)=\sum_{n\geq1}nq^n/(1-q^n)$.
On $\Gm/q^\Z$ the short-model coordinate is
\[
 X(u)=\frac1{12}+
 \sum_{n\in\Z}\frac{q^nu}{(1-q^nu)^2}-2s_1(q).
\]
With $Y=u\,dX/du$ this gives
\[
 Y^2=4X^3-\frac{E_4(q)}{12}X+\frac{E_6(q)}{216},
 \qquad \omega=\frac{dX}{Y}=\frac{du}{u},\quad
 \xi=-X\frac{du}{u}.
\]
Thus $g_2=E_4(q)/12$ and $g_3=-E_6(q)/216$.
These signs agree with complex uniformization by the lattice
$2\pi i\Z+\log(q)\Z$.

\begin{samepage}
\begin{theorem}[Full Tate period matrix]\label{thm:tate-matrix}
For the preceding short model and Tate basis,
\begin{equation}
 P_q=
 \begin{pmatrix}
 t&\ell_q\\
 -E_2(q)t/12&-E_2(q)\ell_q/12-1
 \end{pmatrix}.
 \label{eq:tate-matrix}
\end{equation}
Consequently the class
\begin{equation}
 \nu_q=\xi+\frac{E_2(q)}{12}\omega \label{eq:nu}
\end{equation}
has period row $(0,-1)$.
\end{theorem}
\end{samepage}
\begin{proof}
The pullback of $\omega$ to the torus is $du/u$. The
cyclotomic part of the Tate division tower is the Kummer
tower obtained by taking compatible $p^n$-th roots of $u$.
We recall the $d\log$ comparison on this tower. In the
relative de Rham period sheaf, the compatible system
$u^\flat$ of roots gives the primitive
$\log\bigl(u/[u^\flat]\bigr)$ of $du/u$ for the period
connection. The deck transformation multiplying the roots
by $\varepsilon$ changes this primitive by $-t$.
Subtracting the total \v{C}ech--de Rham differential of the
primitive replaces the form by the negative of this difference.
Its comparison value on $e$ is therefore $t$.
This calculation uses the pro-\'etale Poincar\'e lemma and
the functorial comparison \cite[\S\S6--7]{ScholzeRigid}.
Now \eqref{eq:tate-action} shows that the period on $f$
minus $\ell_q$ is Galois fixed. Both terms lie in $t\BdRp$,
and $K\cap t\BdRp=0$, so the first row is $(t,\ell_q)$.

For the second row consider the convergent rigid theta product
\[
 \Theta(u)=(1-u)\prod_{n\geq1}
 \frac{(1-q^nu)(1-q^n/u)}{(1-q^n)^2},
 \qquad F(u)=u\frac{d}{du}\log\Theta(u).
\]
On annuli avoiding its zeros, termwise differentiation is valid.
Reindexing the product gives $\Theta(qu)=-u^{-1}\Theta(u)$,
so
\[
 F(qu)=F(u)-1,\qquad
 u\frac{dF}{du}
 =-\sum_{n\in\Z}\frac{q^nu}{(1-q^nu)^2}
 =-X(u)+\frac{E_2(q)}{12}.
\]
Thus the meromorphic primitive of $\nu_q$ on the toric
cover is $F$. The constant $F(qu)-F(u)=-1$ describes the
descent of its local primitives under the lattice generator.

To identify its full comparison row, use the semistable
monodromy operator on $\HdR(E_q/K)$.
Coleman--Iovita \cite[Chapters~I--II]{ColemanIovita} identify
this operator with the toric residue followed by the
monodromy pairing and the map from lattice cohomology.
In particular, a class with zero toric residue belongs to
$\ker N$. On an annulus $|q|_p<|u|_p<1$ the primitive $F$
is analytic, so $dF$ has zero annular residue. It follows
that $\nu_q\in\ker N$. The operator $N$ is nonzero because
the toric residue of $du/u$ is $1$ and $v_K(q)>0$.

After extension from the maximal unramified subfield to $K$,
the line $D_{\mathrm{st}}\bigl((V_p(E_q)/\Qp(1))^*\bigr)$
lies in $\ker N$ because the quotient of the Tate representation
is trivial. Since $\ker N$ is one-dimensional, the semistable
comparison identifies $K\nu_q$ with the de Rham
annihilator of the cyclotomic line. Hence
$\langle e,\nu_q\rangle=0$. The basis $(e,f)$ is symplectic
in the convention of~\eqref{eq:legendre-p}, and replacing
$\xi$ by $\nu_q$ preserves the determinant. Therefore
$t\langle f,\nu_q\rangle=-t$ and the row of $\nu_q$ is
$(0,-1)$. Subtracting $E_2(q)\omega/12$ gives~\eqref{eq:tate-matrix}.
\end{proof}

\begin{lemma}\label{lem:E2-unit}
On the entire $p$-adic Tate disk, $E_2(q)$ is a unit and has no zero.
More precisely,
\[
 |E_2(q)-1|_p\leq |24|_p\,|q|_p<1.
\]
\end{lemma}
\begin{proof}
Every $\sigma_1(n)$ is an integer and $|q^n|_p\leq|q|_p$.
The ultrametric inequality gives the bound for the convergent series.
\end{proof}

\subsection{Algebraic independence from the Kummer action}
\begin{theorem}\label{thm:tate-independence}
The elements $t,\ell_q$ are algebraically independent over $\Qpb$.
For the standard de Rham basis $(\omega,\xi)=(du/u,-X\,du/u)$, the two ratios
\[
 \tau_q=\frac{\ell_q}{t},\qquad
 r_q=\frac{-E_2(q)\ell_q/12-1}{-E_2(q)t/12}
\]
are algebraically independent over $\Qpb$.
Both statements remain true after any cycle-basis change in
$\GL_2(\Qp)$, in the corresponding projective coordinates.
\end{theorem}
\begin{proof}
Fix any finite $L/K$.
The cyclotomic character on $G_L$ has infinite image. If $t$ were
algebraic over $L$, a further finite extension would fix $t$, whereas
$\sigma(t)=\chi(\sigma)t$ would force the cyclotomic character to be
trivial there. Hence $t$ is transcendental over $\Qp$.

We claim that $c_q$ is nonzero on $\ker(\chi|_{G_L})$.
Otherwise it factors through the abelian group $\chi(G_L)$.
The cocycle identity and commutativity then give
\[
 (\chi(g)-1)c_q(h)=(\chi(h)-1)c_q(g).
\]
Choosing $g$ with $\chi(g)\ne1$ shows that $c_q(h)$ equals
$(\chi(h)-1)b$ for a fixed $b\in\Qp$. This makes the extension
split over $L$, contradicting Lemma~\ref{lem:kummer}.

Choose $h$ in this kernel with $c=c_q(h)\ne0$.
Then $h^n$ fixes $t$ and sends $\ell_q$ to $\ell_q+nc\,t$ for
every integer $n$. If a nonzero polynomial $A(X,Y)\in L[X,Y]$
vanished at $(t,\ell_q)$, then $A(t,Y)$ would be a nonzero
polynomial, since $t$ is transcendental. It would have the infinitely
many distinct roots $\ell_q+nc\,t$, an impossibility.
Every relation over $\Qpb$ has coefficients in some such finite $L$,
which proves the first assertion.

Put $a=-E_2(q)/12\ne0$. Then
\[
 r_q=\tau_q-\frac1{at},\qquad
 t=\frac1{a(\tau_q-r_q)},\qquad
 \ell_q=\frac{\tau_q}{a(\tau_q-r_q)}.
\]
These birational formulas show that
$\Qpb(t,\ell_q)=\Qpb(\tau_q,r_q)$ and prove the second assertion.
Each change of cycle basis applies the same invertible fractional
linear transformation to the two ratios. Such transformations
preserve their generated field and algebraic independence.
\end{proof}

\begin{corollary}[The unique rational correction]\label{cor:correction}
For $c\in K$, put $\xi_c=\xi+c\omega$.
If $c=E_2(q)/12$, its quasi-period row is $(0,-1)$.
For every other $c$, its quasi-period ratio and the holomorphic period
ratio are algebraically independent over $\Qpb$.
In particular, $\nu_q$ is the unique class of the form $\xi+c\omega$
with rational projective period row.
\end{corollary}
\begin{proof}
The row is $(a_ct,a_c\ell_q-1)$ with $a_c=c-E_2(q)/12$.
If $a_c=0$ the assertion is explicit. Otherwise the birational
argument in Theorem~\ref{thm:tate-independence} applies with $a_c$.
\end{proof}

The standard class $\xi$ never has a $\Qp$-rational projective
period row on the Tate locus, whereas the corrected class $\nu_q$
does. Thus the rationality condition depends on the chosen
second-kind de Rham class.

\subsection{Reduction of the normalized period row to \texorpdfstring{$\Cp$}{Cp}}
Since $\tau_q\in\BdRp$, define
\begin{equation}
 a_q=\theta(\tau_q)\in\Cp . \label{eq:aq}
\end{equation}
\begin{theorem}\label{thm:residue}
The number $a_q$ is transcendental over $\Qp$.
For the standard $\xi$, the associated-graded rows in the Tate
basis are exactly
\[
 \bigl(\theta(W_e/t),\theta(W_f/t)\bigr)=(1,a_q),\qquad
 \bigl(\theta(H_e),\theta(H_f)\bigr)=(0,-1).
\]
After every $\GL_2(\Qp)$ basis change the second projective row is
$\Qp$-rational and the first is transcendental over $\Qp$.
\end{theorem}
\begin{proof}
Equation~\eqref{eq:tate-action} gives
\[
 \sigma(a_q)=\frac{a_q+c_q(\sigma)}{\chi(\sigma)}.
\]
If $a_q$ were algebraic, choose finite $L/K$ fixing it.
Then $c_q(\sigma)=(\chi(\sigma)-1)a_q$ on $G_L$.
Taking a $\sigma$ with $\chi(\sigma)\ne1$ forces $a_q\in\Qp$.
The cocycle is consequently a coboundary, contradicting
Lemma~\ref{lem:kummer}. For the rows, first divide the holomorphic row by $t$ and then
apply $\theta$. Apply $\theta$ to the unscaled second row.
Since $\theta(t)=\theta(\ell_q)=0$,
\eqref{eq:tate-matrix} gives the displayed values. Basis changes apply invertible fractional
linear transformations with $\Qp$ coefficients, which preserve
rationality and transcendence.
\end{proof}

These rows are obtained from the full matrix~\eqref{eq:tate-matrix}
by dividing $W$ by $t$ before applying $\theta$, and by applying
$\theta$ directly to $H$. Thus they record the degree-one holomorphic row, after
trivialization by $t$, and the degree-zero quasi-period row.
The first is the coordinate row of the intrinsic map
\eqref{eq:Fontaine-map}. Proposition~\ref{prop:holomorphic-graded}
identifies its algebraic linear relations with those of the
full holomorphic row. The degree-zero second row records a
different filtration step; reduction alone makes it rational
even though the full standard quasi-period row is not.

\begin{example}[The parameter $q=5$]\label{ex:q5}
Take $K=\Q_5$ and compatible $5$-power roots of $5$.
Then
\[
 t=\log[\varepsilon],\quad
 \ell_5=\log([\widetilde5]/5),\quad
 a_5=\theta(\ell_5/t)
\]
determine $a_5\in\Cp$ once the compatible root towers are fixed.
Theorem~\ref{thm:residue} shows that $a_5$ is transcendental over
$\Q_5$. The integral $q$-expansion gives the separate scalar
calculation
\[
 E_2(5)\equiv 1-24\sum_{n=1}^{7}\sigma_1(n)5^n
       \equiv243581\pmod {5^8}.
\]
Every omitted term is divisible by $5^8$.
This congruence concerns the scalar $E_2(5)$; the root towers
determine the full de Rham periods $t,\ell_5$ and their residue
$a_5=\theta(\ell_5/t)$.
\end{example}

\subsection{Several Kummer parameters and their \texorpdfstring{$1$}{1}-motive}\label{subsec:kummer}
The Tate representation comes from the $1$-motive
$[\Z\to\Gm]$, $1\mapsto q$. For $M=[\Z^r\to\Gm]$,
whose standard generators map to $u_1,\ldots,u_r$, the same
argument computes the full period field. The toric $1$-motive
makes explicit the Kummer extension classes behind
Theorem~\ref{thm:tate-independence}.
Let $u_1,\ldots,u_r\in K^\times$, choose compatible roots, and set
\[
 \ell_i=\log([\widetilde u_i]/u_i)\in t\BdRp.
\]
Let $\kappa_i\in H^1(G_K,\Qp(1))$ be their Kummer classes.

\begin{theorem}[The Kummer rank formula]\label{thm:Kummer-rank}
If $d=\dim_{\Qp}\langle\kappa_1,\ldots,\kappa_r\rangle$, then
\[
 \trdeg_{\Qpb}\Qpb(t,\ell_1,\ldots,\ell_r)=d+1.
\]
In particular, if $\kappa_{i_1},\ldots,\kappa_{i_d}$ form a basis
of the span of the Kummer classes, then
$t,\ell_{i_1},\ldots,\ell_{i_d}$ are algebraically independent
over $\Qpb$.
\end{theorem}
\begin{proof}
First assume the classes are independent.
Restriction to a finite $L/K$ preserves independence: restriction
followed by corestriction is multiplication by $[L:K]$.
Let $c_i$ be the Kummer cocycles on $G_L$, and put
$H=\ker(\chi|_{G_L})$.
The vectors
\[
 \bigl(c_1(h),\ldots,c_r(h)\bigr),\qquad h\in H,
\]
span $\Qp^r$. Otherwise a nonzero linear combination of the
$c_i$ vanishes on $H$. It factors through the abelian group
$\chi(G_L)$ and the cocycle calculation in
Theorem~\ref{thm:tate-independence} makes it a coboundary.
This contradicts independence.

Choose $h_1,\ldots,h_r\in H$ whose vectors form a basis.
Products $h_1^{n_1}\cdots h_r^{n_r}$ fix $t$ and translate
the vector $(\ell_i)$ by $t$ times the corresponding integer
linear combination of these basis vectors.
Suppose a nonzero $A(X,Y_1,\ldots,Y_r)$ over $L$ vanishes at
$(t,\ell_1,\ldots,\ell_r)$.
Transcendence of $t$ makes $A(t,Y_1,\ldots,Y_r)$ nonzero.
After an invertible affine substitution it would vanish at
every point of $\Z^r$.
Induction on the number of variables, using infinitely many
integer roots in one variable, proves that such a polynomial
is zero. This contradiction proves independence over $L$,
and then over $\Qpb$.

For the general case, relabel the classes so that
$\kappa_1,\ldots,\kappa_d$ are independent.
For a dependent class write
$\kappa_j=\sum_{i=1}^d b_i\kappa_i$ with $b_i\in\Qp$.
There is $a\in\Qp$ such that the corresponding cocycles satisfy
$c_j-\sum b_ic_i=(\chi-1)a$.
The element
\[
 z=\ell_j-\sum b_i\ell_i-a t
\]
is $G_K$-fixed, hence lies in $K$.
It lies in $t\BdRp$ as well, so $z=0$ because
$K\cap t\BdRp=0$.
Every remaining $\ell_j$ is thus a $\Qp$-linear combination
of $t,\ell_1,\ldots,\ell_d$.
The upper bound is $d+1$, and the independent case supplies
the matching lower bound.
\end{proof}
The independence required here is $\Qp$-linear independence of
local Kummer classes. Ordinary multiplicative independence over
$\Z$ is not a substitute. Local Kummer theory identifies their
ambient space with the rationalized $p$-adic completion of
$K^\times$, of dimension $[K:\Qp]+1$.
The translation argument identifies a $d$-dimensional unipotent
part of the connected algebraic Galois group and its
one-dimensional cyclotomic quotient. Its dimension $d+1$
agrees with Lemma~\ref{lem:orbit}. Thus the rank of the
local extension classes controls the full period field.

\section{CM and the rational structure on cycles}\label{sec:CM}
\begin{example}[A local CM counterexample]\label{ex:CM}
Let $E:y^2=4x^3-4x$, let $p\equiv1\pmod4$, and choose $i\in\Qp$.
The automorphism $\iota(x,y)=(-x,iy)$ satisfies
\[
 \iota^*\omega=i\omega,\qquad \iota^*\xi=-i\xi.
\]
Write $\iota$ for the induced action on $T_p(E)$.
Since $2,i\in\Z_p^\times$, the operators
\[
 \pi_\pm=\frac12(1\pm i^{-1}\iota)
\]
preserve $T_p(E)$. They are complementary idempotents:
\[
 \pi_\pm^2=\pi_\pm,\qquad
 \pi_+\pi_-=0,\qquad \pi_++\pi_-=1.
\]
Thus
\[
 T_p(E)=T_p(E)^+\oplus T_p(E)^-,\qquad
 T_p(E)^\pm=\ker(\iota\mp i),
\]
The characteristic polynomial of $\iota$ on $T_p(E)$ is $X^2+1$,
so each summand is free of rank one over $\Z_p$.
Choose generators $e_+,e_-$ of these eigenspaces. Endomorphism compatibility
and perfection give
\[
 P=\begin{pmatrix}A&0\\0&B\end{pmatrix},\qquad AB\ne0.
\]
For $\gamma_1=e_++e_-$, $\gamma_2=e_++2e_-$ this becomes
\[
 \begin{pmatrix}A&A\\B&2B\end{pmatrix}.
\]
Thus $H_2/H_1=2$ and $W_2/W_1=1$. Both ratios are rational.
\end{example}
The construction uses the splitting of the CM field at $p$.
The chosen Tate-module basis need not lie in the image of
a prescribed rational Betti homology space.

\begin{proposition}[Betti framing and the CM obstruction]
\label{prop:Betti-CM}
Let $E$ be a CM elliptic curve over a number field.
Fix a complex embedding and a $p$-adic embedding of the number
field, together with extensions to its algebraic closure.
Using these choices, Betti--\'etale comparison identifies
\[
 H_1(E(\C),\Q)\otimes_\Q\Qp \simeq V_p(E).
\]
We regard rational Betti cycles as the image of the untensored
space $H_1(E(\C),\Q)$ under this identification.
A nonzero such cycle cannot annihilate a nonzero de Rham class
defined over a finite $p$-adic field.
In particular, the projective period row of any such class,
computed in a transported rational Betti basis, cannot belong
to $\Pj^1(\Q)$.
\end{proposition}
\begin{proof}
After finite extension all endomorphisms are defined and reduction
is good. Supersingular reduction has no stable line by
Corollary~\ref{cor:irreducible}.
In the ordinary case the only stable lines are the two distinct
CM eigenlines. A rational row relation on a rational Betti vector
$v\ne0$ would make $\Qp v$ one of these lines.

A CM endomorphism outside $\Z$ acts on rational Betti homology
by $M\in M_2(\Q)$ with irreducible quadratic minimal polynomial.
The equality $Mv=\lambda v$ forces $\lambda\in\Q$, since
a nonzero coordinate of $v$ is rational. This contradicts
irreducibility. Comparison respects this endomorphism action.
\end{proof}
This excludes the preceding CM obstruction when the cycle
coordinates come from the chosen rational Betti homology space.
Together with Theorem~\ref{thm:local-A}, the proposition gives
\[
 [H_1:H_2]\in\Pj^1(\Q)
 \quad\Longrightarrow\quad
 [W_1:W_2]\notin\Pj^1(\Qpb)
\]
for elliptic curves over number fields, with rows computed in the
transported Betti basis. This conclusion concerns the local
comparison pairing; transfer of a vanishing complex pairing
requires the additional compatibility discussed in
Section~\ref{subsec:motivic-transfer}.

\section{Ordinary reduction and Katz's Eisenstein series}
\label{sec:ordinary}
Assume $p\geq5$ and good ordinary reduction.
Let $V^+\subset V$ be the connected Tate line introduced in
Corollary~\ref{cor:holomorphic-CM}. The quotient $V/V^+$ is
the rational Tate representation of the \'etale part of the
ordinary $p$-divisible group, so its dual crystalline
realization has slope zero. Exactness and functoriality of
comparison identify
\[
 \Ann_{\HdR(E/K)}(V^+)
   =\DdR\bigl((V/V^+)^*\bigr)=U,
\]
where the annihilator consists of classes whose pairing with
every vector of $V^+$ is zero. Here $U$ is the unit-root line.
For ramified $K$, it is the base extension to $K$ of the
crystalline slope-zero summand.
Our short-model sign convention gives
\begin{equation}
 U=K\left(\xi+\frac{E_{2,p}(E,\omega)}{12}\omega\right).
 \label{eq:unit-root}
\end{equation}
We take~\eqref{eq:unit-root} as the geometric normalization of
$E_{2,p}$, following Katz~\cite[\S5.7, especially Lemma~5.7.8]{Katz}.
With the opposite sign $\eta=-\xi$, the same unit-root formula
and vanishing criterion are \cite[Theorem~2.1 and
Proposition~2.2]{SBp}. The full-row equivalence below combines
this classical formula with Lemma~\ref{lem:rows} and local monodromy.
The Tate $q$-expansion of $E_{2,p}$ is $E_2(q)$, rather than
the $p$-stabilization $E_2(q)-pE_2(q^p)$; this choice fixes the
zero condition studied here. The slope-zero line is the
one-dimensional subisocrystal of Newton slope zero; the
Frobenius iterate that is linear over the unramified base field
has a $p$-adic unit eigenvalue on this line. It is complementary
to the Hodge line $K\omega$.
Mazur--Stein--Tate~\cite[\S3.2]{MST} compute the coefficient
from Frobenius. Over $\Z_p$, for instance, in the basis
$\omega,\eta=x\,dx/y=-\xi$, write
$\varphi^n(\eta)=a_n\omega+b_n\eta$. Then
$E_{2,p}\equiv-12a_n/b_n\pmod {p^n}$.

\begin{theorem}[The ordinary quasi-period condition]
\label{thm:ordinary-condition}
For the standard class $\xi=-x\,dx/y$,
\[
 [H_1:H_2]\in\Pj^1(\Qp)
 \quad\Longleftrightarrow\quad E_{2,p}(E,\omega)=0.
\]
If $E$ is non-CM, this vanishing forces the full holomorphic
period ratio to be transcendental over $\Qp$.
Thus this theorem replaces the complex zero condition in
Proposition~\ref{prop:kumar} by a unit-root zero condition.
\end{theorem}
\begin{proof}
Both conditions are unchanged by finite extension of $K$,
with the comparison bases carried along. In the CM case we
first make such an extension to define the endomorphism action.
For a non-CM ordinary curve, the monodromy description in
Theorem~\ref{thm:dimensions} gives a nonsplit triangular
representation with exactly one stable line, $V^+$. Lemma~\ref{lem:rows} makes
a rational row for $\xi$ equivalent to $K\xi=U$.
For a CM ordinary curve the two annihilator lines are $U$ and
$K\omega$, and $\xi$ cannot span the second.
In either case~\eqref{eq:unit-root} gives the equivalence.
The final assertion is Theorem~\ref{thm:local-A}.
\end{proof}

\begin{corollary}[The holomorphic row throughout the ordinary locus]
\label{cor:ordinary-W}
For a non-CM elliptic curve with good ordinary reduction, the
full holomorphic period ratio is transcendental over $\Qp$,
whether or not its standard quasi-period row is rational.
\end{corollary}
\begin{proof}
Choose a generator $\nu$ of the unit-root line $U$.
It is independent of $\omega$, and its period row is rational
by Lemma~\ref{lem:rows}. Apply Theorem~\ref{thm:local-A}
to the basis $(\omega,\nu)$.
\end{proof}
The vanishing condition in
Theorem~\ref{thm:ordinary-condition} determines when the
\emph{standard} second-kind class has a rational row. The
following complementary statement describes the full pair of
ratios away from the zero locus.

\begin{corollary}[The two ratios away from the ordinary zeros]
\label{cor:ordinary-ratios}
Let $E/K$ have good ordinary reduction without geometric CM
and let $p\ge5$. If $E_{2,p}(E,\omega)\ne0$, the full
holomorphic and standard quasi-period ratios are algebraically
independent over $\Qpb$.
\end{corollary}
\begin{proof}
Put $c=E_{2,p}(E,\omega)/12$ and $\nu=\xi+c\omega$. Choose
a basis $(e_1,e_2)$ of $V$ with $e_1\in V^+$. The
comparison matrix for $(\omega,\nu)$ is
$\left(\begin{smallmatrix}A&C\\0&B\end{smallmatrix}\right)$.
Its three entries $A,C,B$ are
algebraically independent by Theorem~\ref{thm:dimensions}.
The matrix for $(\omega,\xi)$ has second row
$(-cA,B-cC)$, so the two ratios are
\[
 \frac CA,\qquad \frac CA-\frac{B}{cA}.
\]
The triple $(A,C/A,B/A)$ generates the same field as
$(A,C,B)$, and thus has transcendence degree three.
It follows that $C/A$ and $B/A$ are independent, as are
the two displayed ratios. A change of cycle basis acts on
both ratios by invertible fractional linear transformations.
\end{proof}

The correction is compatible with isogenies. For an
ordinary isogeny $f:E\to E'$ of degree $n$, normalize the
differentials by $f^*\omega'=\omega$. In the short models
write $S_f=\sum_{0\ne Q\in\ker f}x(Q)$. The pullback
formula of Sebbar--Besrour
\cite[Proposition~8.7]{SBp} gives
$f^*\xi'=n\xi+S_f\omega$ in our sign convention. In
cycle bases related by $f_*$, the full matrices satisfy
\[
 P_{E'}=
 \begin{pmatrix}1&0\\S_f&n\end{pmatrix}P_E .
\]
The unit-root line pulls back to the unit-root line, so
the coefficients in \eqref{eq:unit-root} satisfy
$c'=nc-S_f$. This expresses the isogeny formula for the
Katz value directly in the full period matrix.

\subsection{The ordinary sigma function and the unit-root correction}
\label{subsec:ordinary-sigma}
Choose an integral short Weierstrass model with unit
discriminant. Use $s=-2x/y$ as formal parameter at the
identity and write $z=\log_E(s)$, so $\omega=dz$ and
$z=s+O(s^2)$. Write $\wp(z)$ for the Laurent expansion
of $x$ in this parameter and put $c=E_{2,p}(E,\omega)/12$.
The Mazur--Tate sigma function, with leading term $z$ after
the change of parameter, has the local expression
\cite{MT,BG}\cite[\S3.1]{MST}
\begin{equation}
 \sigma_p(z)=z\exp\left(
   \int\!\!\int\left(z^{-2}-\wp(z)+c\right)\,dz\,dz\right).
 \label{eq:sigma}
\end{equation}
The two integrations are formal, with zero integration
constants. Since $\wp(z)=z^{-2}+O(z^2)$, the integrand is
regular at $0$ and the expression converges on a sufficiently
small $p$-adic disk.

Let $\zeta_{\mathrm{alg}}(z)=z^{-1}+O(z^3)$ be the formal
primitive with $\zeta_{\mathrm{alg}}'=-\wp(z)$.
Differentiation of~\eqref{eq:sigma} gives
$\zeta_p=\sigma_p'/\sigma_p=\zeta_{\mathrm{alg}}+cz$ and
\begin{equation}
 d\zeta_p=\left(-\wp(z)+\frac{E_{2,p}(E,\omega)}{12}\right)dz
   =\xi+\frac{E_{2,p}(E,\omega)}{12}\omega .
 \label{eq:sigma-correction}
\end{equation}
The right side is a global second-kind differential whose
class spans the unit-root line. Changing $c$ to $c+a$
multiplies the local sigma expression by $\exp(az^2/2)$ and
adds $a\omega$ to the second-kind class. The unit-root
splitting fixes the coefficient that gives a rational row
in Theorem~\ref{thm:ordinary-condition}. On the Tate locus the
corresponding correction is $E_2(q)/12$ by
Theorem~\ref{thm:tate-matrix}.

\subsection{Multiplicity of ordinary zeros}
Work on a fine tame level $N\geq3$, prime to $p$.
Multiplicity means order on the smooth modular curve, or on
its ordinary Igusa cover, with respect to any local trivialization
of the Hodge bundle. It is not order in the coarse $j$-coordinate
at its ramification points $0,1728$.

Put $P=E_{2,p}$, $Q=E_4$, and $R=E_6$ in this subsection.
On the ordinary Igusa cover, a trivialization of the connected
$p$-divisible group supplies a canonical invariant differential
and hence trivializes the Hodge bundle. In this trivialization,
Katz's theta operator is a local derivation $D$
\cite[Lemma~5.8.1]{Katz}.
If $q_{\mathrm{ST}}\in1+\mathfrak m$ is a Serre--Tate coordinate
on an ordinary deformation disk, then $D$ is, up to a unit,
$q_{\mathrm{ST}}\,d/dq_{\mathrm{ST}}$. It is nonvanishing
because $q_{\mathrm{ST}}$ is a unit. This is why nonzero
values of $DP$ and $D^2P$ detect the first and second
orders of vanishing used below: in a local parameter $s$,
$D=v(s)\,d/ds$ with $v(0)\ne0$, and $2$ is invertible
because $p\ge5$.
Its $q$-expansion is $q\,d/dq$. Kodaira--Spencer and the
$q$-expansion principle~\cite[\S\S5.7--5.8]{Katz} give
\begin{equation}
 DP=\frac{P^2-Q}{12},\qquad
 DQ=\frac{PQ-R}{3},\qquad
 DR=\frac{PR-Q^2}{2}.
 \label{eq:Ramanujan}
\end{equation}
These identities hold over the characteristic-zero ordinary locus
and after reduction modulo $p$.

We recall the multiplicity classification of Sebbar--Besrour
\cite[Proposition~3.1, Theorem~3.3, and Corollary~3.8]{SBp}.
\begin{theorem}[Ordinary zero multiplicities, Sebbar--Besrour]\label{thm:zeros}
Every ordinary zero of $E_{2,p}$ is simple unless $j=0$.
Every ordinary point with $j=0$ is a zero of multiplicity two.
Such points occur precisely for $p\equiv1\pmod3$.
The same multiplicities hold on the characteristic-$p$ ordinary locus.
\end{theorem}
\begin{proof}
These are \cite[Proposition~3.1 and Corollary~3.8]{SBp} in
characteristic zero and \cite[Theorem~3.3]{SBp} after reduction.
Their argument uses the nonvanishing Serre--Tate derivation and
the Ramanujan identities~\eqref{eq:Ramanujan}; the extra
automorphisms give the zeros at the ordinary elliptic points.
The multiplicity convention above is the fine-level convention
used in those results.
\end{proof}
An ordinary $j=1728$ point is likewise a zero, by the order-four
automorphism. It is simple since $Q\ne0$. It is ordinary exactly
when $p\equiv1\pmod4$.

\subsection{The CM zeros}
The following lemma concerns the almost holomorphic form
$E_2^*(\tau)=E_2(\tau)-3/(\pi\operatorname{Im}\tau)$, rather than
the holomorphic $E_2$ in Question~\eqref{eq:A}.

\begin{lemma}[Sebbar--Besrour]\label{lem:E2star}
The zeros of $E_2^*$ are exactly the $\mathrm{SL}_2(\Z)$-orbits
of $i$ and $\rho=-1/2+i\sqrt3/2$.
\end{lemma}
\begin{proof}
This is \cite[Lemma~8.1]{SBp}. Its proof uses the weight-two
transformation law to reduce to the standard fundamental domain,
then Fourier estimates to exclude zeros away from the two elliptic
orbits. We use the classification without repeating those estimates.
\end{proof}

\begin{proposition}[CM zero criterion, Sebbar--Besrour]\label{prop:CM-zero}
For an elliptic curve with geometric CM and good ordinary reduction
at $p\geq5$,
\[
 E_{2,p}(E,\omega)=0\quad\Longleftrightarrow\quad
 j(E)\in\{0,1728\}.
\]
\end{proposition}
\begin{proof}
The result is \cite[Theorem~8.2 and Corollary~8.3]{SBp}.
A geometric CM curve over a $p$-adic field has algebraic
$j$-invariant and is homothetic over $\Cp$ to a CM model over
$\Qb$. The vanishing condition is unchanged by that homothety,
so the cited algebraic CM classification applies.
\end{proof}
The proof in \cite{SBp} compares the two CM eigenspaces defined
over $\Qb$: the complementary eigenspace becomes antiholomorphic
under complex comparison, and the Legendre relation reduces
vanishing to Lemma~\ref{lem:E2star}. This comparison uses
algebraically defined eigenspaces. Transport of the two arbitrary
marked relations of Question A is considered separately in
Subsection~\ref{subsec:motivic-transfer}.

\subsection{A non-CM zero over \texorpdfstring{$\Q_{11}$}{Q11}}
\begin{example}[The ordinary Katz zero over $\Q_{11}$]
\label{ex:eleven}
There is a unique $j_{11}\in7+11\Z_{11}$ for which
\begin{equation}
 E_j:\quad y^2=4x^3-\frac{3j}{j-1728}x-\frac{j}{j-1728},
 \qquad \omega_j=dx/y
 \label{eq:Ej}
\end{equation}
has $E_{2,11}(E_j,\omega_j)=0$. It is the only ordinary
Katz-zero $j$-invariant in $\C_{11}$.
These assertions specialize Sebbar--Besrour's classification
\cite[Theorem~3.7 and Corollaries~3.8 and~3.10]{SBp}.
The curve $E_{j_{11}}/\Q_{11}$ has good ordinary reduction
without geometric CM. Its standard full quasi-period row is
$\Q_{11}$-rational and its full holomorphic period ratio is
transcendental over $\Q_{11}$.

For $p=11$, the supersingular polynomial and the two elliptic
exponents of \cite[\S1 and Theorem~3.7]{SBp} are
\[
 S_{11}(X)=X(X-1),\qquad\delta_{11}=\varepsilon_{11}=1.
\]
Their polynomial for nonelliptic ordinary residue zeros becomes
\[
 Z_{11}(X)=-12(2X-1)+8(X-1)+6X=X-7
 \quad\text{in }\F_{11}[X].
\]
Theorem~3.7 and Corollary~3.10 of \cite{SBp} give its unique
ordinary lift $j_{11}\in7+11\Z_{11}$, over $\Q_{11}$ itself.
Their count is $\lfloor11/12\rfloor+1=1$
\cite[Corollary~3.8]{SBp}, and neither elliptic
class is ordinary at $11$, so this is the only ordinary zero
invariant. For comparison with our short-model normalization,
evaluation with $u=j/(j-1728)$ gives $Q=36u$, $R=-216u$,
and $A=E_{10}(E_j,\omega_j)=QR$. Here $E_{10}=E_4E_6$ is
a weight-ten form, and its reduction is the Hasse invariant.
On the residue disk $j\equiv7$, one has $u\equiv3$ and
$A\equiv9\pmod{11}$. The normalized weight-twelve form then gives
\begin{equation}
 E_{2,11}\equiv E_{12}/A
 =\frac{5(j-7)}{j-1}\pmod{11}.
 \label{eq:11-congruence}
\end{equation}
Thus this model has the same residual zero as the cited
classification. Its discriminant is
$46656j^2/(j-1728)^3$, a unit throughout the disk, and its
Hasse value is a unit, so it has good ordinary reduction.
Since $j_{11}$ is neither $0$ nor $1728$, Proposition~\ref{prop:CM-zero}
excludes geometric CM. The full quasi-period row is rational
by Theorem~\ref{thm:ordinary-condition}, and the holomorphic
period ratio is transcendental by Theorem~\ref{thm:local-A}.
\end{example}
The equation $E_{2,11}(E_j,\omega_j)=0$ together with $j\in7+11\Z_{11}$ uniquely determines $j_{11}$.
The Frobenius algorithm of~\cite{MST}, followed by Hensel iteration,
computes further digits.
In a basis with $H_1=0$, the quasi-period ratio is $\infty$.
Replacing $(e_1,e_2)$ by $(e_2,e_1+r e_2)$ realizes any prescribed
$r\in\Q_{11}$ as a finite quasi-period ratio, while the period
ratio stays transcendental.

For a general prime $p\ge5$, Sebbar--Besrour
\cite[Theorems~3.5 and~3.7, Corollary~3.10]{SBp} classify the
ordinary zeros and lift each simple nonelliptic residue zero
uniquely, including over ramified residue disks. Their zeros have
exactly $\lfloor p/12\rfloor+1$ distinct $j$-invariants, all in
finite unramified extensions of $\Qp$; the unramified degree of
a nonelliptic lift equals the degree of its residual invariant.
Together with Theorem~\ref{thm:ordinary-condition}, this classifies
the ordinary $j$-invariants for which the standard full
quasi-period row is rational. Theorem~\ref{thm:local-A} gives
the transcendental holomorphic ratio at every nonelliptic zero.
The full period field gives a stronger joint conclusion.

\begin{corollary}[Full and reduced ordinary period fields]
\label{cor:ordinary-zero-periods}
Let $E/K$ have good ordinary reduction without geometric CM,
and suppose $p\ge5$. Use the unit-root basis
$(\omega,\nu)$ with $\nu=\xi+E_{2,p}(E,\omega)\omega/12$
and a cycle basis $(e_1,e_2)$ with $e_1\in V^+$.
Let $h=\langle e_2,\nu\rangle$ and $\rho=W_2/W_1$.
Then $h$ and $W_1/t$ are units of $\BdRp$ and
$\rho\in\BdRp$. The fields of the full and normalized matrices are
\begin{equation}
 \Qpb(P)=\Qpb(t,h,\rho),\qquad
 \Qpb(P_{\mathrm{gr}})=\Qpb\bigl(\theta(h),\theta(\rho)\bigr).
 \label{eq:ordinary-fields}
\end{equation}
The three generators in the first field and the two generators
in the second are algebraically independent over $\Qpb$.
At a Katz zero one has $\nu=\xi$, so both assertions apply to
the standard matrix. In particular they apply to
Example~\ref{ex:eleven}.
\end{corollary}
\begin{proof}
Since $U$ annihilates $V^+$, the lower-left entry of $P$
vanishes. Write $\delta=\det P/t\in K^\times$. Then
$W_1h=\delta t$. Both $W_1/t$ and $h$ belong to $\BdRp$,
and their product is a unit, so both are units. The formulas
\[
 P=\begin{pmatrix}\delta t/h&\delta t\rho/h\\0&h\end{pmatrix},
 \qquad
 P_{\mathrm{gr}}=
 \begin{pmatrix}
   \delta/\theta(h)&\delta\theta(\rho)/\theta(h)\\
   0&\theta(h)
 \end{pmatrix}
\]
give the field inclusions in \eqref{eq:ordinary-fields}.
Conversely $t=\det P/\delta$, $h$ is the lower-right entry
of $P$, and $\rho=W_2/W_1$. The lower-right entry of
$P_{\mathrm{gr}}$ is $\theta(h)$, and the ratio of its
two top entries is $\theta(\rho)$.
The ordinary Borel cases of Theorems~\ref{thm:dimensions}
and~\ref{thm:graded-relation-ideal} give the respective
transcendence degrees three and two.
\end{proof}

\section{Characteristic \texorpdfstring{$p$}{p}: crystalline and Drinfeld realizations}
\label{sec:drinfeld}
\subsection{What crystalline comparison does provide}\label{subsec:crystalline}
Let $E_0$ be an elliptic curve over a finite field $\kappa$ of
characteristic $p$. Put $K_0=W(\kappa)[1/p]$, where $W(\kappa)$
is the Witt-vector ring, and set
\[
 C(E_0)=H^1_{\mathrm{cris}}(E_0/W(\kappa))[1/p].
\]
This is a two-dimensional $K_0$-space with semilinear Frobenius.
Choose a finite extension $K/K_0$, with residue field
$\kappa_K$, and a smooth proper elliptic curve
$\widetilde E/\mathcal O_K$ whose special fiber is identified
with $E_0\otimes_\kappa\kappa_K$. Write $E/K$ for its
generic fiber.
Crystalline--de Rham comparison gives
\begin{equation}
 C(E_0)\otimes_{K_0}K \simeq \HdR(E/K).
 \label{eq:cris-lift}
\end{equation}
This comparison is indispensable for the unit-root construction
in~\eqref{eq:unit-root}; see \cite[Theorem V.2.3.2]{Berthelot}
and \cite[Section 9.1]{BC}.

It does not, by itself, supply the analogue of the two complex
cycle directions. The $K$-linear comparison~\eqref{eq:cris-lift} is an
isomorphism of $K$-vector spaces, so its matrix in $K_0$ and $K$
bases has entries in $K$. These entries are algebraic over $\Qp$.
A transcendence question over $\Qp$ for that matrix cannot be the
full local question~\eqref{eq:Ap}. Moreover, the crystalline
space does not carry a distinguished rank-two $\Q$-cycle lattice.
Its Frobenius eigendirections are additional local structure,
rather than a replacement for rational Betti homology.

The Hodge line is part of the choice of lift. Different lifts of
the same ordinary $E_0$ share its isocrystal while their Hodge
filtrations vary. If one includes Frobenius and the filtration and
then compares with the Tate representation through $\Bcris$,
one returns to the mixed-characteristic period theory of
Section~\ref{sec:local}. When the lift $E$ is defined over
$\Qb$, one can retain its de Rham basis over $\Qb$ and compare
the crystalline actions of rational algebraic correspondences on
$E_0$ with this global structure. This gives a transcendence
problem over $\Q$. Section~\ref{subsec:crystalline-consequence} develops
a concrete consequence in that setting.

The geometric \'etale $p$-adic Tate module
$\varprojlim_n E_0[p^n](\overline\kappa)$ has rank one when
$E_0$ is ordinary and rank zero when it is supersingular. Thus it
cannot replace the rank-two $V_p(E)$ of Section~\ref{sec:local}.
The characteristic-$p$ assertion in Theorem~\ref{thm:zeros}
concerns zeros on the ordinary modular curve; it is not an
elliptic rank-two period implication over a finite field.

\subsection{Why Drinfeld modules enter}
Fix a prime power $q=p^m$ and put
\[
 A=\F_q[\theta],\quad k=\F_q(\theta),\quad
 k_\infty=\F_q((1/\theta)),\quad
 \C_\infty=\widehat{\overline{k_\infty}}.
\]
Let $\overline k$ be the algebraic closure of $k$ inside
$\C_\infty$.

Drinfeld modules are function-field analogues of elliptic curves
for the present purpose: a rank-two module has an entire
exponential, a discrete rank-two $A$-lattice, and
a two-dimensional de Rham realization with quasi-periods.
They are not elliptic curves over a finite field, and their
quasi-periodic functions are entire rather than meromorphic zeta
functions. Their lattice and de Rham realization give a
two-by-two comparison matrix whose holomorphic and quasi-period
rows can be studied separately.

We work in generic $A$-characteristic, meaning that the structural
map from $A$ into the coefficient field is injective. The coefficient field
still has characteristic $p$. This use of generic
characteristic must be distinguished from finite
$A$-characteristic, where the structural map has a nonzero kernel.

\subsection{Lattices, quasi-periods, and algebraic coefficients}
Write $\mathcal T$ for the $q$-power operator,
so $\mathcal T c=c^q\mathcal T$.
Consider a rank-two Drinfeld module over $\overline k$,
\[
 \phi_\theta=\theta+g\mathcal T+\Delta\mathcal T^2,\qquad
 g,\Delta\in\overline k,\quad \Delta\ne0.
\]
Its exponential, normalized by $\exp_\phi'(0)=1$, is entire and has kernel
$\Lambda=A\omega_1+A\omega_2$.
The biderivation $\delta_\theta=\mathcal T$ gives the entire
quasi-periodic function $F_\delta$, characterized by
\[
 F_\delta(\theta z)-\theta F_\delta(z)=\exp_\phi(z)^q,\qquad
 F_\delta(0)=F'_\delta(0)=0.
\]
Its restriction to $\Lambda$ is $A$-linear. Put
$\eta_i=F_\delta(\omega_i)$ and
\[
 P_\phi=\begin{pmatrix}\omega_1&\omega_2\\
                         \eta_1&\eta_2\end{pmatrix}.
\]
The inner class $\delta^{(1)}_a=a-\phi_a$ has function
$z-\exp_\phi(z)$, whose value at $\omega_i$ is $\omega_i$
because $\exp_\phi(\omega_i)=0$. It therefore supplies the
holomorphic period row.
The classes $\delta^{(1)}$ and $\delta$ form a basis of
the rank-two de Rham realization; see
\cite[Section~3.1]{CP}. Gekeler's de Rham comparison
\cite{Gekeler} therefore makes $P_\phi$ invertible. Any second de Rham class defined over
$\overline k$ and independent of the inner class $\delta^{(1)}$
can be used in place of $\delta$.

The rational cycle space is $k\Lambda=k\otimes_A\Lambda$,
of dimension two over $k$; it retains the lattice's natural
$k$-structure. Let $F_\phi=\End(\phi)\otimes_A k$ and
$s=[F_\phi:k]$.
For rank two, $s=1$ or $2$. We say that $\phi$ has CM when
$s=2$, allowing inseparable quadratic endomorphism fields.
The period theorem of Chang--Papanikolas~\cite[Theorem 1.2.2]{CP}
gives
\begin{equation}
 \trdeg_{\overline k}\overline k(P_\phi)=4/s.
 \label{eq:CP}
\end{equation}
The theorem, with Conrad's appendix, includes inseparable
endomorphism fields. In characteristic two a quadratic
endomorphism field can be purely inseparable, so the CM
argument below does not assume diagonalizability.

\begin{proposition}[The Drinfeld zero reformulation]
\label{prop:drinfeld-reformulation}
For a rank-two Drinfeld lattice $\Lambda$ over $\C_\infty$,
use the standard biderivation $\delta_\theta=\mathcal T$.
Its quasi-period row is in $\Pj^1(k)$ if and only if a
primitive lattice period has zero quasi-period.
Equivalently, after changing an $A$-basis and scaling the
first period to $1$, the resulting lattice $A+Az$ satisfies
$\mathcal E(z)=0$, where
\[
 \mathcal E(z)=F_{\delta,z}(1),\qquad
 z\in\Omega=\C_\infty\setminus k_\infty,
\]
and $F_{\delta,z}$ is the normalized quasi-periodic function of
the Drinfeld module attached to $A+Az$.
The change of period ratio is fractional linear over $k$,
so it preserves algebraicity over $k$.
\end{proposition}
\begin{proof}
Clear the denominators of a nonzero $k$-linear relation between
the quasi-periods and remove the common divisor in
$A=\F_q[\theta]$. The resulting coprime pair $(a,b)$ extends
to a matrix in $\GL_2(A)$, since $A$ is a principal ideal
domain. The primitive period $a\omega_1+b\omega_2$ has
zero quasi-period. Conversely, such a period gives a
nonzero $k$-linear relation.

For a homothety $\Lambda\mapsto c\Lambda$, the exponential
satisfies $\exp_{c\Lambda}(cw)=c\exp_\Lambda(w)$.
The defining functional equation and its normalization then give
$F_{\delta,c\Lambda}(cw)=c^qF_{\delta,\Lambda}(w)$.
Scaling the chosen primitive period to $1$ therefore
preserves its vanishing. Finally, an invertible fractional
linear transformation over $k$ has an inverse over $k$.
Algebraicity of the ratio is preserved in both directions.
\end{proof}
This is the change-of-basis counterpart of
Proposition~\ref{prop:kumar}. The notation $\mathcal E$ records
a quasi-period function; it does not identify that function with
the classical $E_2$ or Katz's $E_{2,p}$.
By Gekeler's formula \cite[Theorem~7.10]{GekelerQuasi},
there is a constant $c_{\mathrm G}\in\C_\infty^\times$,
independent of $z$, such that
\[
 \mathcal E(z)=c_{\mathrm G}\,
      \frac{E_{\mathrm{Dr}}(z)}{h(z)}.
\]
Here $E_{\mathrm{Dr}}$ is the Drinfeld false Eisenstein
series and $h$ is the cusp form of weight $q+1$ with no
zeros on $\Omega$. The constant records the normalization
of these functions; its nonvanishing and that of $h$ give
the precise equivalence
$\mathcal E(z)=0\Longleftrightarrow E_{\mathrm{Dr}}(z)=0$.
The definition $\mathcal E(z)=F_{\delta,z}(1)$ fixes the
quasi-period normalization used in the proposition.
Under the homothety $\Lambda\mapsto c\Lambda$, the module
coefficients change by
\[
 (g,\Delta)\longmapsto
 \bigl(c^{1-q}g,c^{1-q^2}\Delta\bigr).
\]
Normalizing a lattice period to $1$ can therefore take coefficients
in $\overline k$ to coefficients outside $\overline k$.
The algebraic-coefficient hypothesis must be checked for the
normalized module when applying~\eqref{eq:CP}.

\subsection{A linear row theorem from Yu's quasi-period theorem}
\label{subsec:Yu-rows}
Yu's sub-$t$-module theorem \cite{Yu97} is the function-field
analytic subgroup theorem relevant to linear period relations.
For the row condition needed here, his earlier quasi-period
theorem already suffices and gives an assertion in arbitrary
rank.

For clarity, the algebraic de Rham group
$H_{\mathrm{DR}}(\phi/\overline k)$ is the quotient of
algebraic biderivations
\[
 \delta:A\longrightarrow\overline k\{\mathcal T\}\mathcal T,
 \qquad \delta_{ab}=a\delta_b+\delta_a\phi_b,
\]
by the strictly inner biderivations
$\delta^{(m)}_a=am-m\phi_a$, where
$m\in\overline k\{\mathcal T\}\mathcal T$.
The braces denote the ring of additive polynomials, with the
commutation rule $\mathcal Tc=c^q\mathcal T$.
Strictly inner biderivations have zero quasi-periods on the
lattice, so the period functional depends only on the de Rham
class. The inner class $\delta^{(1)}$ used above remains
nonzero in this quotient; the constant polynomial $1$ is
excluded from the strictly inner subspace. See
\cite[Section~3]{Yu90} and \cite[Section~3.1]{CP}.

\begin{proposition}[Yu's theorem and linear independence of a row]
\label{prop:Yu-row}
Let $\phi$ be a Drinfeld $A$-module of rank $r\ge1$ in generic
$A$-characteristic over $\overline k$. Let
$0\ne\delta\in H_{\mathrm{DR}}(\phi/\overline k)$, and let
$\omega_1,\ldots,\omega_r$ be an $A$-basis of its period
lattice. Then
\[
 F_\delta(\omega_1),\ldots,F_\delta(\omega_r)
 \quad\text{are linearly independent over }k.
\]
This assertion has no restriction on the endomorphism field.
\end{proposition}
\begin{proof}
Yu's quasi-period theorem \cite[Theorem~3.1]{Yu90} states that
$F_\delta(u)$ is transcendental over $k$ whenever $u\ne0$,
$\exp_\phi(u)\in\overline k$, and $\delta$ has nonzero de
Rham class. In particular $F_\delta(\lambda)\ne0$ for every
nonzero lattice period $\lambda$.

Suppose $\sum_i b_iF_\delta(\omega_i)=0$ with $b_i\in k$,
not all zero. Clear denominators to obtain $a_i\in A$,
not all zero, satisfying the same relation.
The period $\lambda=\sum_i a_i\omega_i$ is nonzero because
the $\omega_i$ form an $A$-basis. The $A$-linearity of
$F_\delta$ on the lattice gives $F_\delta(\lambda)=0$,
contradicting Yu's theorem.
\end{proof}

For rank two, the proposition excludes a $k$-rational period
row for every nonzero algebraic de Rham class. This includes
the standard quasi-period class and every second class used
below. The matrix argument that follows records the
endomorphism constraints explicitly and gives the stronger
non-CM algebraic-independence conclusion from \cite{CP}.

\begin{lemma}[Endomorphisms in de Rham coordinates]\label{lem:Drin-endo}
Suppose $s=2$, and choose $u\in\End(\phi)$ whose image in
$F_\phi$ does not belong to $k$.
Its action on $k\Lambda$ is $M\in M_2(k)$ with irreducible
quadratic minimal polynomial. There are
$\lambda,\mu,b\in\overline k$ with $\lambda,\mu\notin k$ such that
\begin{equation}
 P_\phi M=
 \begin{pmatrix}\lambda&0\\b&\mu\end{pmatrix}P_\phi .
 \label{eq:Drin-endo}
\end{equation}
The two diagonal entries may coincide in the inseparable case.
\end{lemma}
\begin{proof}
The action on $k\Lambda$ is the regular two-dimensional action
of $k(u)$, so its minimal polynomial is irreducible quadratic.
The functorial de Rham comparison intertwines this action with
the action on the algebraic de Rham group.
The line of inner classes is invariant: its period row sends a
lattice element to itself, and applying $u$ multiplies this row
by the scalar differential of $u$, denoted $\lambda$.
The matrix in a basis starting with that inner class is therefore
lower triangular over $\overline k$. Its characteristic polynomial
equals that of $M$, because comparison is invertible.
Consequently its other diagonal entry $\mu$ is the second root
of the same irreducible polynomial, and neither root belongs to $k$.
This uses only the naturality and algebraic definition of the
de Rham group; see~\cite{Gekeler,CP}.
\end{proof}

\begin{theorem}[Quasi-period rows over algebraic Drinfeld coefficients]
\label{thm:drinfeld-A}
For a rank-two Drinfeld module over $\overline k$, let the
second de Rham class be defined over $\overline k$ and independent
of the inner class. Its projective quasi-period row lies outside
$\Pj^1(k)$. If the module has no CM, the two ratios
$\omega_2/\omega_1$ and $\eta_2/\eta_1$ are algebraically
independent over $\overline k$.
\end{theorem}
\begin{proof}
Proposition~\ref{prop:Yu-row} gives the rational-row exclusion
for all endomorphism fields. In the non-CM case, the four
matrix entries are algebraically independent by
\cite[Theorem~1.2.2]{CP}. Clearing denominators in a polynomial
relation of the ratios would give a nonzero polynomial relation
of the four entries.
\end{proof}

In the CM case, the endomorphism calculation identifies the
algebraic projective period points.

\begin{proposition}[Projective period points in the CM case]
\label{prop:Drin-CM-rows}
Use the notation of Lemma~\ref{lem:Drin-endo} and write
$W=(\omega_1,\omega_2)$ and $H=(\eta_1,\eta_2)$.
The holomorphic point belongs to
$\Pj^1(k(\lambda))\setminus\Pj^1(k)$.
The quasi-period point is algebraic over $k$ precisely when the
second de Rham class is an eigendirection for $u$, that is, when
$b=0$ in~\eqref{eq:Drin-endo}. In that case it belongs to
$\Pj^1(k(\mu))\setminus\Pj^1(k)$.
If the endomorphism field is inseparable, the quasi-period point
is transcendental over $k$ for every second class independent
of the inner class.
\end{proposition}
\begin{proof}
The first row of~\eqref{eq:Drin-endo} gives $WM=\lambda W$.
The left eigenspace has dimension one and is defined over
$k(\lambda)$. A $k$-rational eigenvector would force
$\lambda\in k$. This proves the first assertion. If $b=0$,
the same argument applies to $H$ and $\mu$.

Suppose $b\ne0$ and $H=h v$ for some $h\ne0$ and
$v\in\overline k^2\setminus\{0\}$. The second row gives
\begin{equation}
 h(vM-\mu v)=bW .
 \label{eq:CM-drinfeld}
\end{equation}
All four matrix entries then lie in $\overline k\,h$,
so their field has transcendence degree at most one. This
contradicts the value two in~\eqref{eq:CP}. Finally, when the
endomorphism field is inseparable one has $\lambda=\mu$.
The matrix in~\eqref{eq:Drin-endo} cannot be scalar because
it is similar to $M$, whose minimal polynomial has degree two.
Thus $b\ne0$ in every such basis.
\end{proof}

\begin{corollary}[The Drinfeld $j$-invariant at a quasi-period zero]\label{cor:drinfeld-j}
For a rank-two Drinfeld module over $\C_\infty$, using the
standard biderivation $\delta_\theta=\mathcal T$, a
$k$-rational quasi-period row forces its invariant
\[
 j_{\mathrm{Dr}}(\phi)=g^{q+1}/\Delta
\]
to be transcendental over $k$. In particular,
$\mathcal E(z)=0$, equivalently $E_{\mathrm{Dr}}(z)=0$,
forces $j_{\mathrm{Dr}}(z)\notin\overline k$.
\end{corollary}
\begin{proof}
If $j_{\mathrm{Dr}}\in\overline k$, the module is
isomorphic over $\C_\infty$ to one over $\overline k$.
Indeed, if $g\ne0$, conjugate by a scalar $c$ with
$c^{q-1}=g$ to obtain coefficients $1$ and
$1/j_{\mathrm{Dr}}$. If $g=0$, use
$c^{q^2-1}=\Delta$ to obtain coefficients $0$ and $1$.
Under the corresponding lattice homothety the two
standard quasi-periods are multiplied by $c^q$, as in
Proposition~\ref{prop:drinfeld-reformulation}.
Their projective $k$-rationality is preserved,
contradicting Theorem~\ref{thm:drinfeld-A}.
\end{proof}
This is the Drinfeld counterpart of
Proposition~\ref{prop:zero-j}. It concerns the moduli invariant
$j_{\mathrm{Dr}}(z)$. For the lattice parameter itself the
following degree exclusion is unconditional.

\begin{corollary}[Quadratic parameters are excluded]
\label{cor:drinfeld-quadratic}
An algebraic zero $z\in\Omega$ of $\mathcal E$ has
$[k(z):k]\ge3$.
\end{corollary}
\begin{proof}
Since $z\in\Omega$, it lies outside $k$. If $[k(z):k]=2$,
choose $a\in A\setminus\{0\}$ clearing the denominators of the
coefficients expressing $z^2$ in the basis $(1,z)$ over $k$.
Then multiplication by $az$ preserves $A+Az$ and gives an
endomorphism of its Drinfeld module outside the image of $A$. The lattice
has CM, so its singular modulus is algebraic by
\cite[Theorem~3]{GekelerSingular}. The algebraicity statement
in that theorem covers inseparable quadratic fields.
This contradicts Corollary~\ref{cor:drinfeld-j}.
\end{proof}
Whether $\mathcal E$ can vanish at an algebraic parameter
of higher degree remains open.

\begin{example}[CM with its natural lattice]
Take $\phi_\theta=\theta+\mathcal T^2$.
It is the rank-one Carlitz module for $\F_{q^2}[\theta]$,
viewed as rank two over $A$. Its lattice is
$\F_{q^2}[\theta]\omega$.
For $\zeta\in\F_{q^2}\setminus\F_q$ take
$\omega_1=\omega$, $\omega_2=\zeta\omega$.
Since the exponential is $\F_{q^2}$-linear, uniqueness in the
quasi-periodic functional equation gives
$F_\delta(\zeta z)=\zeta^qF_\delta(z)$.
Thus
\[
 \frac{\omega_2}{\omega_1}=\zeta,\qquad
 \frac{\eta_2}{\eta_1}=\zeta^q\notin k.
\]
The nonzero determinant is
$(\zeta^q-\zeta)\omega\eta_1$.
This is consistent with transcendence degree two and with
Theorem~\ref{thm:drinfeld-A}. The natural $k$-structure on
the lattice prevents the split-CM obstruction obtained by allowing
arbitrary $\Qp$-bases in Example~\ref{ex:CM}.
\end{example}

\section{Analytic zeta and sigma functions and the de Rham class}
\label{sec:analytic}
\subsection{The Tate primitive}
In this section, a $p$-adic zeta function means an elliptic
Weierstrass zeta-type primitive, rather than a $p$-adic
interpolation of the Riemann zeta function.
A formal elliptic logarithm is an analytic coordinate on a
sufficiently small identity neighborhood. It does not provide a
rank-two additive period lattice in $\Cp$.

For a Tate curve, the uniformization is $\Gm/q^\Z$.
Let $F(u)=u\,d\log\Theta(u)/du$ be the theta derivative used
in Theorem~\ref{thm:tate-matrix}. Choose a branch of $\log_p$ on $K^\times$ and its
compatible extension to algebraic extensions. On local analytic
domains of the toric covering avoiding the zeros of $\Theta$,
put
\[
 Z_q(u)=F(u)-\frac{E_2(q)}{12}\log_p(u).
\]
The calculation in that theorem gives
\[
 dZ_q=\xi,\qquad
 Z_q(qu)-Z_q(u)=-1-\frac{E_2(q)}{12}\log_p(q).
\]
The function $Z_q$ is a local analytic primitive with the
chosen logarithm branch. Its increment depends on that branch
and is one scalar analytic increment, whereas~\eqref{eq:tate-matrix} is a
two-column comparison matrix. For the customary branch
$\log_p(p)=0$, the scalar $\log_p(p)$ vanishes but
$\ell_p/t$ is transcendental by
Theorem~\ref{thm:tate-independence}. The two constructions
record different data.

\subsection{The ordinary zeta primitive}
The ordinary sigma calculation in
Subsection~\ref{subsec:ordinary-sigma} gives
$d\zeta_p=\xi+E_{2,p}(E,\omega)\omega/12$ in the formal
logarithm coordinate. The Tate theta primitive belongs to
the toric covering, and the ordinary primitive is defined on
a formal neighborhood of the identity. In both cases the
differential represents the class whose full comparison row
annihilates the connected Tate direction.
The logarithms of algebraic points in the next subsection
use the analytic group logarithm on the universal vector
extension.

\subsection{Fuchs--Pham and logarithms on the universal vector extension}
\label{subsec:analytic-subgroups}
The preceding primitives are related to algebraic de Rham
classes. Analytic subgroup theorems address a complementary
arithmetic question: which algebraic linear constraints can the
logarithm of an algebraic point satisfy? We fix an embedding
$\Qb\hookrightarrow\Cp$ throughout this subsection.

For a commutative algebraic group $G/\Qb$, let $G(\Cp)_f$
be the domain of its canonical $p$-adic logarithm
\[
 \log_G:G(\Cp)_f\longrightarrow\Lie(G)\otimes_{\Qb}\Cp.
\]
A point belongs to this domain when some strictly increasing
sequence of positive integer multiples tends to the identity.
Fuchs--Pham \cite[Theorem~2.2]{FP15} give the following form
of the $p$-adic analytic subgroup theorem, restated in
\cite[Theorem~1.1]{FP16}.

\begin{theorem}[The $p$-adic analytic subgroup theorem]
\label{thm:Fuchs-Pham}
Let $G/\Qb$ be a commutative algebraic group of positive
dimension, and let $L\subseteq\Lie(G)$ be a nonzero
$\Qb$-linear subspace. If
\[
 P\in G(\Qb)\cap G(\Cp)_f,\qquad
 0\ne\log_G(P)\in L\otimes_{\Qb}\Cp,
\]
then there is an algebraic subgroup $H\subseteq G$ of positive
dimension, defined over $\Qb$, such that
$\Lie(H)\subseteq L$ and $P\in H(\Qb)$.
\end{theorem}

Fuchs--Pham credit Bertrand's earlier formulation and Matev's
proof, and provide an independent proof. Their sequel \cite{FP16} applies the theorem to
elliptic and abelian logarithms. The period pairings with
specified rational structures in \cite{MS} are compared in
Subsection~\ref{subsec:motivic-transfer}.
The universal vector extension
in \eqref{eq:universal-extension} gives the following
consequence, with the same subgroup obstruction as in
Proposition~\ref{prop:complex-column}.

\begin{proposition}[A logarithm ratio on the universal vector extension]
\label{prop:analytic-log-slope}
Let $E/\Qb$, and fix algebraic coordinates $(X,Y)$ on
$\Lie(E^\natural)$ with $d\pi$ equal to the first coordinate.
For
\[
 P\in E^\natural(\Qb)\cap E^\natural(\Cp)_f,
 \qquad \log_{E^\natural}(P)=(u,v),
\]
if $u\ne0$, then $u,v$ are linearly independent over $\Qb$.
Equivalently, $v/u$ is transcendental over $\Q$.
\end{proposition}
\begin{proof}
A nonzero relation $au+bv=0$ over $\Qb$ has $b\ne0$.
The line $L=\ker(aX+bY)$ projects nontrivially to
$\Lie(E)$. Theorem~\ref{thm:Fuchs-Pham} produces a
positive-dimensional subgroup with Lie algebra contained in
$L$. Its identity component has dimension one and Lie algebra
$L$, so it projects nontrivially to $E$.
Lemma~\ref{lem:universal-subgroups} excludes this subgroup.
\end{proof}

The coordinate $v$ is a logarithm coordinate on an algebraic
extension of $E$ by $\Ga$. Expressing it through a Weierstrass
zeta primitive requires an algebraic lift of the elliptic point
and a compatible trivialization of the additive fibre.
These choices specify the comparison with the analytic zeta
values studied in \cite{Pham}. The unit-root correction in
\eqref{eq:sigma-correction} specifies a de Rham class over
the local field $K$.
Proposition~\ref{prop:ordinary-arithmetic} addresses the
arithmetic of that coefficient over $\Q$ under the stated
crystalline period prediction.

The fields of coefficients separate this proposition from
Theorem~\ref{thm:local-A}. In algebraic coordinates, the
$p$-adic logarithm of an algebraic point belongs to a finite
extension of $\Qp$ containing the point and the coefficients
of the group. Proposition~\ref{prop:analytic-log-slope} asserts
transcendence over $\Q$. Theorem~\ref{thm:local-A} concerns
full comparison periods and transcendence over $\Qp$.
Its rational vanishing relation defines a Galois-stable
Tate-module line by Lemma~\ref{lem:rows}; the line can come
from a connected $p$-divisible group, with no corresponding
elliptic algebraic subgroup. These are the respective
geometric mechanisms of the analytic and comparison arguments.

\begin{example}[Global independence and local Kummer rank]
\label{ex:logarithm-fields}
For $p=5$, set
\[
 \beta=\frac{\log_5(81)}{\log_5(16)}\in\Q_5.
\]
Both arguments lie in $1+5\Z_5$, where the logarithm is
injective, so both logarithms are nonzero.
The number $\beta$ is transcendental over $\Q$.
Indeed, if $\beta\in\Qb$, apply
Theorem~\ref{thm:Fuchs-Pham} to $(16,81)\in\Gm^2(\Q)$
and its logarithm line. The resulting subgroup has a
one-dimensional torus as identity component. A nonzero
integer character trivial on that torus gives
\[
 a\log_5(16)+b\log_5(81)=0,
 \qquad (a,b)\in\Z^2\setminus\{(0,0)\}.
\]
Injectivity on $1+5\Z_5$ gives $16^a81^b=1$.
Unique factorization then forces $a=b=0$, a contradiction.

The two local Kummer classes nevertheless span a
one-dimensional $\Q_5$-space: they belong to the
rationalized principal-unit group of $\Q_5$, on which the
logarithm is an isomorphism to $\Q_5$.
Theorem~\ref{thm:Kummer-rank} therefore gives transcendence
degree two for their full comparison periods together with $t$.
This example distinguishes the global algebraic coefficients in
the analytic subgroup theorem from the local coefficient field
in the Kummer rank formula.
\end{example}

\section{The return to the complex problem}\label{sec:bridge}
\subsection{The generalized period conjecture implies Question A}
\label{subsec:motivic}
For an elliptic curve in a short Weierstrass model over a
finitely generated field $F\subset\C$, the relevant conjectural
bound is~\eqref{eq:motivic-intro}. Replacing $F$ by an algebraic extension does not change
$\trdeg_\Q F(P)$, where $P$ is the period matrix.
For elliptic curves the motivic group has the same dimension
as the Mumford--Tate group, namely $4$ without CM and $2$
with CM. This is the elliptic case of Andr\'e's generalized
period conjecture, explained in
\cite{Andre,Bertolin02,Bertolin20}.
The field of coefficients is essential to its formulation.

Sebbar \cite[Proposition~13.1]{Sebbar} proves the conditional
implication for a rational projective quasi-period row.
We recall the argument using the normalized lattice.
\begin{proposition}[Conditional consequence of Grothendieck--Andr\'e]
\label{prop:motivic-A}
Assume~\eqref{eq:motivic-intro} for non-CM elliptic curves over
finitely generated subfields of $\C$. Then Question A holds.
\end{proposition}
\begin{proof}
By Proposition~\ref{prop:kumar}, a counterexample can be
written as a lattice $\Z+\Z\alpha$ with $\alpha\in\Qb\cap\HH$
and $E_2(\alpha)=0$. Its period and quasi-period rows are
$(1,\alpha)$ and $(0,-2\pi i)$.
Proposition~\ref{prop:zero-j} excludes CM.
For its Weierstrass model take
$F=\Q(\alpha,g_2,g_3)$. Its standard differentials are
$F$-defined and
\[
 F(\omega_1,\omega_2,\eta_1,\eta_2)
   =\Q(\alpha,g_2,g_3,2\pi i).
\]
Since $\alpha$ is algebraic, the right-hand side has
transcendence degree at most $3$ over $\Q$.
The conjectural lower bound is $4$, a contradiction.
\end{proof}
The argument uses the generalized bound over the coefficient
field of the normalized lattice. It gives the conditional
consequence directly; the simultaneous local compatibility
considered below is a separate possible approach.

\subsection{A geometric form of the rational relation}
For the universal vector extension
\eqref{eq:universal-extension} over $\C$, use the de Rham
coordinates introduced in Section~\ref{subsec:complex-subgroups}.
Its analytic uniformization is
\[
 E^{\natural,\mathrm{an}}\simeq
 \C^2/\bigl(\Z(\omega_1,\eta_1)+\Z(\omega_2,\eta_2)\bigr),
\]
with projection to $E$ in the first coordinate.
Let $\mathfrak h=\{(z,0):z\in\C\}$.

\begin{proposition}[A Zariski-dense analytic subgroup]
\label{prop:leaf}
If $[\eta_1:\eta_2]\in\Pj^1(\Q)$, then the analytic subgroup
$\exp_{E^\natural}(\mathfrak h)$ is isomorphic as a complex
Lie group to $\C^\times$ and is Zariski dense in $E^\natural$.
\end{proposition}
\begin{proof}
The kernel of the additive quasi-period homomorphism
$\Lambda\to\C$ has rank one over $\Z$.
Its rank cannot be two by the Legendre relation.
Consequently $\mathfrak h$ meets the exponential lattice in
$\Z(\omega_0,0)$ for some $\omega_0\ne0$.
The analytic subgroup is $\C/\omega_0\Z\simeq\C^\times$.

The Zariski closure of the image is a connected algebraic
subgroup of $E^\natural$. Since $\mathfrak h$ projects
nontrivially to $\Lie(E)$, its closure is not contained in
the additive kernel. Lemma~\ref{lem:universal-subgroups}
therefore makes that closure all of $E^\natural$.
\end{proof}
The complex rational quasi-period relation thus produces a
one-dimensional, nonalgebraic analytic subgroup.
This supplies the geometric form of the rational relation.
The same algebraic subgroup obstruction appears in
Propositions~\ref{prop:complex-column} and~\ref{prop:analytic-log-slope}.
Sebbar~\cite[Section~14 and Proposition~14.1]{Sebbar} considers
a related homomorphism into $\Ga^2\times E$ under algebraicity
of both period ratios. Its Zariski closure is $\Ga\times E$,
and the numerical rank hypothesis of Waldschmidt's theorem fails.
Here the ambient group is the nontrivial universal vector extension
$E^\natural$, and the rational quasi-period relation gives a
rank-one exponential kernel. The two constructions exhibit
different algebraic subgroups and different obstructions.
Proposition~\ref{prop:zero-j} places our unrestricted curve over
a field with transcendental $j$, so an arithmetic algebraicity
criterion requires further input in that setting.

Bost's theorem~\cite[Theorem~2.3]{Bost} gives a precise possible type of input:
for an algebraic group over a number field, a rational Lie
subalgebra is algebraic if its reductions are closed under the
restricted $p$-power operation at almost every prime.
For a Lie algebra of vector fields in characteristic $p$,
this operation is the $p$-fold iterate of a derivation.
Translation-invariant complex leaves have universal cover
$\C$ in this rank-one case. Bounded-above subharmonic
functions on that cover are constant, giving the Liouville
condition required by Bost's criterion.
We obtain the following immediate obstruction.

\begin{corollary}\label{cor:Bost}
Let $E$ and a line $\mathfrak h\subset\Lie(E^\natural)$
be defined over a number field $F$. Choose a smooth integral
model of $E^\natural$ and extend $\mathfrak h$ as a line
subbundle of its Lie algebra after excluding finitely many
primes of $F$. If $\mathfrak h$ projects nontrivially to
$\Lie E$, its reductions fail to be closed under restricted
$p$-th powers at infinitely many remaining good primes.
\end{corollary}
\begin{proof}
Otherwise Bost's criterion would make $\mathfrak h$ the Lie
algebra of an algebraic subgroup. Lemma~\ref{lem:universal-subgroups} excludes that subgroup.
\end{proof}
For a hypothetical complex counterexample,
Proposition~\ref{prop:zero-j} gives transcendental $j$.
Its curve therefore does not descend to a number field.
An analogous approach in that case would require a model over a
finitely generated arithmetic base, an extension of the specified
Lie line, and an algebraicity criterion for families with the
corresponding restricted-power hypotheses verified.

\subsection{The two relations that must be transported}
Suppose a complex counterexample to~\eqref{eq:A} existed.
Adjoin its algebraic ratio $\alpha=\omega_2/\omega_1$
to a finitely generated field $F\subset\C$ containing $g_2,g_3$.
The counterexample would give
\begin{equation}
 \eta_2-r\eta_1=0,\qquad
 \omega_2-\alpha\omega_1=0,\qquad r\in\Q,
 \label{eq:two-relations}
\end{equation}
after a rational cycle-basis change if necessary.
This includes the projective case by choosing a suitable chart.
Proposition~\ref{prop:zero-j}, applied after
Proposition~\ref{prop:kumar}, shows that its elliptic curve
is geometrically non-CM.

Choose an embedding $F\hookrightarrow K$ into a finite
$p$-adic field. Such embeddings exist: send a transcendence
basis to an algebraically independent tuple in $\Qp$ and
extend to the finite algebraic extension into $\Qpb$. The
finitely many images of its algebraic generators lie in a
single finite extension $K$ of $\Qp$. With compatible
embeddings of the algebraic closures, Betti--\'etale comparison gives
\[
 H_1(E(\C),\Q)\otimes_\Q\Qp \simeq V_p(E_K),
\]
where $E_K$ is obtained from the model over $F$ through
$F\hookrightarrow K$. This transports the rational Betti cycles
and preserves their rational linear combinations. Vanishing
of a pairing involves the de Rham class as well, so the
equalities that would be needed are
\begin{equation}
 H_2-rH_1=0,\qquad W_2-\alpha_pW_1=0.
 \label{eq:needed-transfer}
\end{equation}
Here $\alpha_p$ is the image of $\alpha$ in $K$, which is
algebraic over $\Qp$. The equalities concern the local
comparison of the same marked cycles and classes.

\begin{proposition}[A sufficient compatibility condition in $B_2$]
\label{prop:transfer}
For a hypothetical complex counterexample, it is sufficient
to find a single local realization with the transported cycles
and de Rham classes above for which
\begin{equation}
 q_2(H_2)-r\,q_2(H_1)=0,\qquad
 q_2(W_2)-\alpha_p\,q_2(W_1)=0
 \quad\text{in }B_2.
 \label{eq:truncated-transfer}
\end{equation}
For these transported Tate vectors and de Rham classes, the equations are
equivalent to the full relations \eqref{eq:needed-transfer}.
Existence of such a realization for every hypothetical
counterexample would prove Question A.
\end{proposition}
\begin{proof}
Proposition~\ref{prop:linear-truncation} lifts
\eqref{eq:truncated-transfer} to
\eqref{eq:needed-transfer}. The base-changed curve remains
geometrically non-CM: geometric endomorphisms are unchanged
by extension of algebraically closed characteristic-zero
fields. The first relation gives a $\Qp$-rational
quasi-period row and the second an algebraic holomorphic row.
Comparison keeps the matrix invertible. These conditions
contradict Theorem~\ref{thm:local-A}.
\end{proof}

The two equations in~\eqref{eq:truncated-transfer} are the
required extra input. Transport of the cycle vectors alone
does not give compatibility of even a single homogeneous
period relation, as the following example shows.

\begin{example}[Failure of transport of a holomorphic relation]
\label{ex:failed-transport}
Take the lattice $\Z+\Z\alpha$ with $\alpha=i\sqrt[3]{2}$
and its Weierstrass model over $F=\Q(\alpha,g_2,g_3)$.
Since $-\alpha^2=2^{2/3}$, the algebraic number $\alpha$ has
degree greater than two over $\Q$. Schneider's theorem
\cite{Schneider} makes $j(\alpha)$ transcendental over $\Q$.
Include $j(\alpha)$ in a transcendence basis of $F$ and
choose its image in $p^{-1}\Z_p^\times$ to be transcendental.
Choose the other images successively in $\Qp$ to be
algebraically independent over the preceding ones. These
choices are possible because the fields generated at each
step have countable algebraic closures and the indicated
$p$-adic open sets are uncountable. Extend the embedding to
$F$ with image in a finite extension $K/\Qp$.

The resulting curve has potentially multiplicative reduction,
since $v_p(j(E_K))<0$. In the transported cycle basis
Corollary~\ref{cor:holomorphic-CM} gives
$[W_1:W_2]\notin\Pj^1(\Qpb)$. In particular,
$W_2-\alpha_pW_1\ne0$, although $\omega_2-\alpha\omega_1=0$
in the complex comparison.
\end{example}

The local groups in Table~\ref{tab:local} give a related
difference in the relation ideals. For instance, the
connected Tate direction of an ordinary non-CM curve yields
a linear relation and leaves a full local period field of
transcendence degree three over $\Qpb$. The expected complex
degree over algebraic coefficients is four.

\subsection{Motivic relations and simultaneous transfer}
\label{subsec:motivic-transfer}
The two vanishing relations can be described as equations on a
comparison matrix. Fix the field $F$, the two chosen Betti cycles,
and the de Rham basis used above. For a generic matrix $X$ put
\[
 R=F[X_{11},X_{12},X_{21},X_{22},\det(X)^{-1}],
\]
and define $f_H=X_{22}-rX_{21}$ and
$f_W=X_{12}-\alpha X_{11}$.
Complex and local comparisons give evaluation maps
$R\to\C$ and $R\to\BdR$, using the respective embeddings
of $F$. Simultaneous transfer asks whether both specified
elements of the first kernel belong to the second kernel.
For the test modulo $t^2$ in \eqref{eq:truncated-transfer}, use
the polynomial subalgebra $F[X_{11},X_{12},X_{21},X_{22}]$:
$q_2(P)$ has noninvertible determinant, so it is not a point
of $\GL_2(B_2)$. The equations $f_H,f_W$ belong to that
polynomial subalgebra, and
Proposition~\ref{prop:linear-truncation} detects their full
vanishing there.

There is a precise sufficient motivic mechanism. Suppose the
elliptic motive, specified cycle vectors, and de Rham classes
admit compatible tensor realizations in a category of motives. Its period torsor is the
space of tensor isomorphisms between the two fiber functors.
An identity in the coordinate ring of that torsor is respected
by every comparison point, including a compatible local point.
Thus, if the equations $f_H=0$ and $f_W=0$ hold in the coordinate
ring of that period torsor, both relations transfer. The substantive issue is establishing that they are
such identities from their vanishing at one complex point.

For a non-CM elliptic motive alone the group is $\GL_2$.
A function $X\mapsto \lambda(Xv)$, with
$\lambda\ne0$ a de Rham covector and $v\ne0$ a cycle,
is not identically zero on $\GL_2$: an invertible matrix can
send $v$ to a vector outside $\ker\lambda$. A zero pairing
at one comparison point therefore does not, on its own, define
a tensor identity. Marked $1$-motives or relative arithmetic
data would have to supply additional structure.

Over algebraic coefficient fields, the strong period conjecture
expresses genericity of the complex comparison point in its
period torsor, and is a mechanism for recognizing all its
algebraic relations. For a coefficient field not algebraic over $\Q$,
the generalized bound~\eqref{eq:motivic-intro} is a statement
about total transcendence degree over $\Q$. It does not
assert genericity relative to every enlarged field $F$.
For example, normalizing a lattice to have first period $1$
creates a relation over its coefficient field without producing
a new elliptic endomorphism.

Consequently the standard motivic conjectural picture has two
distinct uses here. Its generalized complex bound already
excludes a counterexample by Proposition~\ref{prop:motivic-A}.
A separate approach through Proposition~\ref{prop:transfer}
would need a realization-compatible explanation for both
marked relations. Conjectures about Hodge or Tate classes
concern algebraic tensors; applying them requires identifying
the relevant tensors and justifying specialization.

The period spaces of \cite{MS} offer a concrete comparison of
the domains in such a construction. Their $B_2$-valued pairings
are defined using specified $\Q$-structures and restricted
de Rham argument spaces over number fields with good reduction
\cite[Definitions~4.7--4.11]{MS}. Their formal period relations
are generated by bilinearity and functoriality of $1$-motives.
Presentation depth two, in \cite[Definitions~4.5--4.6]{MS}, refers
to exact sequences involving at most two copies of a motive;
it is distinct from reduction modulo $t^2$ in $B_2$.
Our comparison below concerns the defined spaces of arguments.

For actual Tate vectors and de Rham classes,
Proposition~\ref{prop:linear-truncation} proves independently that
the two local equations are determined in $B_2$, for every
reduction type. The following intrinsic facts about $\Phi_E$
then specify which constant Tate vectors can belong to the
logarithm-based spaces.

First assume only that $E/K$ has good reduction, and use
the map $\Phi_E$ of \eqref{eq:Fontaine-map}.
Fix a basis $e$ of the one-dimensional $\Qp$-space
$\Qp(1)$, and write $\Qpb(1)=\Qpb\,e\subset\Cp(1)$.
This specifies the subspace
$\Lie(E)\otimes_K\Qpb(1)$ of the target of $\Phi_E$.
We use the same basis $e$ to place formal logarithms in this
twisted target below.

\begin{proposition}[Algebraic coefficients of Fontaine periods]
\label{prop:Fontaine-algebraic}
For an elliptic curve with good reduction,
\[
 \Phi_E(V_p(E))\cap
 \bigl(\Lie(E)\otimes_K\Qpb(1)\bigr)=\{0\}.
\]
\end{proposition}
\begin{proof}
Choose a $K$-basis of $\Lie(E)$ and suppose that
$\Phi_E(v)=ae\ne0$ with $a\in\Qpb$.
Choose finite $L/K$ containing $a$, and let $W$ be the
$\Qp$-span of the $G_L$-orbit of $v$.
Equivariance gives
\[
 \Phi_E(\sigma v)=\chi(\sigma)ae
 \qquad(\sigma\in G_L).
\]
Thus $\Phi_E(W)=\Qp\,ae\simeq\Qp(1)$ is a subquotient
of $V_p(E)|_{G_L}$.

Good reduction makes $V_p(E)$ crystalline, and the crystalline
category is closed under subquotients. Duality would give
$\Qp(-1)$ as a subquotient of
$H^1_{\mathrm{et}}(E_{\overline L},\Qp)$.
Let the residue field of $L$ have order $q=p^f$.
The linear Frobenius iterate on crystalline $H^1$ would have
eigenvalue $q$. Its characteristic polynomial is
$T^2-a_LT+q$, where
$a_L=q+1-\#E_0(\F_q)$. If one root is $q$, the other is
$1$, giving $\#E_0(\F_q)=0$. The identity point makes this
impossible.
\end{proof}

Let $\mathcal E/\mathcal O_K$ be the smooth proper model,
let $\mathcal G^0$ be the connected component of its
$p$-divisible group, and put
\[
 V^0=T_p(\mathcal G^0)\otimes_{\Z_p}\Qp.
\]

\begin{lemma}\label{lem:Fontaine-connected}
One has $\ker\Phi_E=\Ann_{V_p(E)}(\omega)$ for any nonzero
invariant differential $\omega$. The restriction of $\Phi_E$
to $V^0$ is injective.
\end{lemma}
\begin{proof}
The first assertion follows from
Proposition~\ref{prop:holomorphic-graded}, or directly from
Lemma~\ref{lem:finite-filtration}: vanishing of the
degree-one pairing means
$\langle v,\omega\rangle\in t^2\BdRp$, hence full vanishing.
For supersingular reduction $V^0=V_p(E)$ and
Corollary~\ref{cor:irreducible} excludes a nonzero rational
kernel. For ordinary reduction $V^0=V^+$ has rank one, and
its de Rham annihilator is the unit-root line, complementary
to $K\omega$, as in \eqref{eq:unit-root}.
It has no nonzero vector annihilated by $\omega$.
\end{proof}

Now let $E$ be defined over a number field $\mathbb K$, with
good reduction at the chosen completion $K$.
The space $h_p(E)$ in \cite[Definition~4.7]{MS} is the
$\Q$-span of formal-group logarithms of globally algebraic
points. Each point and the group are defined over a finite
extension of the local field, and the convergent logarithm
takes its values in that complete field. We identify the logarithm of a point $P$ with
$\log_E(P)\otimes e$ using the fixed basis of $\Qp(1)$.
With this convention,
\begin{equation}
 h_p(E)\subset\Lie(E)\otimes_K\Qpb(1).
 \label{eq:MS-local-algebraicity}
\end{equation}
This is algebraicity over $\Qp$; the values need not be
algebraic over $\Q$.

By \cite[Definition~4.8]{MS}, the space
$\mathcal H_p^\varphi(E)$ is the inverse image of $h_p(E)$
under the scalar extension of $\Phi_E$ to
$V_p(E)\otimes_{\Qp}\Cp$.
Their refined structure
$\widetilde{\mathcal H}_p^\varpi(E)$ is contained in
$V^0\otimes_{\Qp}B_2$, and its reduction to
$V_p(E)\otimes_{\Qp}\Cp$ lies in
$\mathcal H_p^\varphi(E)$
\cite[Definitions~4.9--4.10]{MS}. Only this containment and
the reduction property enter the following proposition; it
applies to any argument subspace with those two properties.

\begin{proposition}[Specified structures and constant Tate vectors]
\label{prop:MS-intersection}
Use the constant embeddings $v\mapsto v\otimes1$ in the
indicated scalar extensions. Then
\[
 \mathcal H_p^\varphi(E)\cap V_p(E)=\ker\Phi_E,
 \qquad
 \widetilde{\mathcal H}_p^\varpi(E)\cap V_p(E)=\{0\}.
\]
The second intersection is taken inside
$V_p(E)\otimes_{\Qp}B_2$.
\end{proposition}
\begin{proof}
Membership of $v$ in the first intersection says
$\Phi_E(v)\in h_p(E)$.
Equation~\eqref{eq:MS-local-algebraicity} and
Proposition~\ref{prop:Fontaine-algebraic} make this equivalent
to $\Phi_E(v)=0$. The converse uses $0\in h_p(E)$.

A constant vector in the refined structure belongs to $V^0$:
its image in $(V_p(E)/V^0)\otimes B_2$ is zero, and scalar
extension from $\Qp$ is faithful. Its reduction belongs to
$\mathcal H_p^\varphi(E)$, so the first equality gives
$\Phi_E(v)=0$. Lemma~\ref{lem:Fontaine-connected} then
gives $v=0$.
\end{proof}

The refined cycle space in
Proposition~\ref{prop:MS-intersection} has trivial intersection
with the constant Tate vectors in both ordinary and supersingular
reduction. This distinction concerns the actual arguments of
the pairing, even in the supersingular case, where $V^0=V_p(E)$.

The de Rham argument of the refined pairing
\cite[Definition~4.11]{MS} lies in the nonzero-slope part
$\widetilde N(E)$ of crystalline cohomology, with coefficient
field $\mathbb K^u=K^{\mathrm{ur}}\cap\overline{\mathbb K}$.
Their Fontaine pairing uses the specified $\Qb$-defined
holomorphic classes. At a Katz zero, \eqref{eq:unit-root}
puts the standard second-kind class $\xi$ on the slope-zero
line, outside $\widetilde N(E)$.

For Example~\ref{ex:eleven}, algebraicity of $j_{11}$ would
supply a number-field model. The standard class would still
span the slope-zero line, and the transported cycles would
still be constant Tate vectors. A construction that transfers
the two relations must respect both arguments of each
pairing and accommodate the finitely generated coefficient
field of a hypothetical complex counterexample. The next
subsection treats global algebraicity of $j_{11}$ through
crystalline endomorphism periods.
\subsection{A crystalline prediction for the ordinary example}
\label{subsec:crystalline-consequence}
Crystalline comparison becomes an arithmetic period theory when
one fixes a de Rham basis over algebraic numbers. We can use it
to sharpen the question about $j_{11}$.

Let $E/\Qb$ have good ordinary reduction $E_0$ over
$\F_{p^f}$, and fix an embedding $\Qb\hookrightarrow\Qpb$.
Choose a finite extension $K/\Qp$ over which the chosen model
and its good reduction are defined.
Let $\pi$ be the $p^f$-power Frobenius endomorphism of $E_0$.
Through~\eqref{eq:cris-lift}, its pullback is a $K$-linear
endomorphism of $\HdR(E/K)$. Write $\Phi_\pi$ for its
matrix in a $\Qb$-defined basis, using columns for the images
of basis vectors. This is the action of the $p^f$-power
Frobenius endomorphism; the underlying crystalline
$p$-Frobenius is semilinear. Its entries lie in $K\subset\Qpb$ and need not lie in
the embedded $\Qb$.
Its characteristic polynomial is
$T^2-aT+p^f$ with $a\in\Z$; its two eigenvalues
$\alpha,\beta\in\Qb$ are distinct, and the unit eigenvalue
cuts out the unit-root line.

For an ordinary elliptic curve over a finite field, the
geometric rational endomorphism algebra is the imaginary
quadratic field
\[
 \End(E_{0,\overline{\F}_p})\otimes_\Z\Q=\Q(\pi).
\]
In particular every rational endomorphism is $a+b\pi$ with
$a,b\in\Q$. Its matrix in crystalline--de Rham comparison
is $aI+b\Phi_\pi$. Hence the entries of $\Phi_\pi$ generate
the field of endomorphism periods of the reduction in the
chosen algebraic de Rham basis. This identifies the matrix
prediction below with the ordinary endomorphism-motive case
of \cite[Examples~5.14 and~6.17]{AF}.

\begin{conjecture}[Ordinary crystalline period prediction]
\label{conj:ordinary-crystalline}
If $E/\Qb$ is non-CM and $E_0$ is ordinary, then
\begin{equation}
 \trdeg_{\Qb}\Qb(\Phi_\pi)=2.
 \label{eq:ordinary-crystalline}
\end{equation}
The base $\Qb$ here denotes its specified image in $\Qpb$.
\end{conjecture}
This is the ordinary elliptic endomorphism case of the
crystalline period conjectures of Ancona--Fr\u{a}\c{t}il\u{a}.
For $M=\End(h^1(E))$, their de Rham group is
$\operatorname{PGL}_2$, the ordinary crystalline group is a
one-dimensional torus $T$, and the period space is
$H_M\simeq\operatorname{PGL}_2/T$, of dimension two
\cite[Examples~5.14 and~6.17]{AF}. After extending the algebraic
coefficient field, it parametrizes ordered pairs of distinct
Frobenius eigenlines. The transported Frobenius matrix, whose
two eigenvalues are algebraic and distinct, determines that pair
and is determined by it.

Sebbar--Besrour~\cite[Remark~8.11 and Proposition~8.12]{SBp}
use \cite[Conjecture~8.3 and Proposition~8.4]{AF}, which assert
genericity of the Berthelot point in this period space, to
exclude the unit-root line $K\xi$. Our matrix formulation has
the same conjectural content in this case. Indeed,
\cite[Example~6.17]{AF} gives
$P_p^{\mathrm{mot}}(M)=\mathcal O(H_M)$, so their
Conjecture~8.1 for the single motive and Conjecture~8.3 for its
generated tensor category coincide here. The required
unramifiedness condition is automatic by \cite[Remark~7.7]{AF}.
This condition concerns CM objects in the tensor category;
it differs from the unramified fields of definition of the
zeros in \cite[Corollary~3.10]{SBp}. The base in
\eqref{eq:ordinary-crystalline} is the embedded $\Qb$;
the local full-period degree $d_p(E)=3$ in
Table~\ref{tab:local} has base $\Qpb$.

The same period-space mechanism excludes every algebraic
Katz value, extending the conditional nonvanishing of
\cite[Proposition~8.12]{SBp}. The next statement records this
strengthening and the joint information in both eigenlines.
\begin{proposition}[Conditional arithmetic of the ordinary zeros]
\label{prop:ordinary-arithmetic}
Assume Conjecture~\ref{conj:ordinary-crystalline}.
For a non-CM elliptic curve over $\Qb$ with good ordinary
reduction at $p\ge5$, choose a short Weierstrass model and differential
$\omega$ over $\Qb$. Then the value
$E_{2,p}(E,\omega)$ is transcendental over $\Q$.
More precisely, both Frobenius eigenlines admit generators
$\xi+c_0\omega$ and $\xi+c_1\omega$, with $c_0=E_{2,p}/12$
on the unit-root line, and $c_0,c_1$ are algebraically
independent over the embedded $\Qb$.
In particular every nonelliptic ordinary Katz-zero invariant,
including $j_{11}$ of Example~\ref{ex:eleven}, is
transcendental over $\Q$ under this conjecture.
\end{proposition}
\begin{proof}
Write $F$ for the embedded $\Qb$. Let
$D=\Dcris(V_p(E)^*)$ over $K_0$, the maximal unramified
subfield of $K$, and let $D=D_0\oplus D_1$ be its slope
decomposition. The slope lines after extension to $K$
are the eigenlines $L_0,L_1$ of $\Phi_\pi$, with
$L_0$ the unit-root line. Equation~\eqref{eq:unit-root}
gives $L_0=K(\xi+c_0\omega)$ with $c_0=E_{2,p}/12$.

The line $L_1$ differs from the Hodge line already
unconditionally. Indeed, if $L_1=K\omega$, the projection
$D\to D_1$ along $D_0$ would commute with crystalline
Frobenius and preserve the Hodge filtration on $D\otimes_{K_0} K$.
By full faithfulness of $\Dcris$
\cite[Proposition~9.1.11]{BC}, this rank-one idempotent would
split $V_p(E)^*$ into two one-dimensional representations.
Serre's ordinary criterion \cite[Chapter~IV, Appendix~A.2.4]{Serre}
would then give geometric CM. Thus $L_1=K(\xi+c_1\omega)$
for a unique $c_1$, and $c_0\ne c_1$.

With $\alpha$ the unit eigenvalue, put
\[
 S=\begin{pmatrix}c_0&c_1\\1&1\end{pmatrix},\qquad
 \Phi_\pi=S\begin{pmatrix}\alpha&0\\0&\beta\end{pmatrix}S^{-1}.
\]
Hence $F(\Phi_\pi)\subseteq F(c_0,c_1)$. Conversely, the
spectral projectors
$(\Phi_\pi-\beta I)/(\alpha-\beta)$ and
$(\Phi_\pi-\alpha I)/(\beta-\alpha)$ recover the two
eigenlines over $F(\Phi_\pi)$. Taking the ratio of the
coordinates of any nonzero column gives $c_0,c_1$.
Consequently $F(\Phi_\pi)=F(c_0,c_1)$, and the conjectural
transcendence degree two proves their algebraic independence.
In particular $E_{2,p}=12c_0$ is transcendental over $\Q$.

If a nonelliptic zero invariant were algebraic over $\Q$,
it would admit a model and differential over $\Qb$, still with
good ordinary reduction after finite extension. The CM
classification of Proposition~\ref{prop:CM-zero} makes that
curve non-CM. Its Katz value is zero, contradicting the first
assertion. This applies to $j_{11}$.
\end{proof}
Geometrically, prescribing any algebraic value $c_0=a$ fixes
the first line in the ordered pair parametrized by
$\operatorname{PGL}_2/T$, a proper divisor. This is the
divisor argument of \cite[Remark~8.11 and
Proposition~8.12]{SBp}, with an arbitrary algebraic line in
place of $K\xi$. The matrix calculation makes explicit the
joint coordinates of the generic point. Unconditionally, the
arithmetic of the non-CM zero invariants over $\Q$ remains open.

\subsection{Concrete next questions}
The results leave several focused research questions.
\begin{question}
Is the invariant $j_{11}$ in Example~\ref{ex:eleven}
algebraic or transcendental over $\Q$?
The curve $E_{j_{11}}$ has no geometric CM. This is the
$p=11$ instance of the global algebraicity question in
\cite[Question~8.13]{SBp}.
Proposition~\ref{prop:ordinary-arithmetic} gives transcendence
conditionally on the ordinary crystalline period prediction.
\end{question}
\begin{question}
For a finitely generated coefficient field supporting a hypothetical
complex zero of $E_2$ at algebraic $\tau$, which properties of
the annihilated rational cycle constrain the unit-root or crystalline
lines of its reductions? In particular, can one obtain the
restricted-power hypotheses needed for an algebraicity criterion?
\end{question}
\begin{question}
Can a realization-compatible construction on a universal
vector extension or a $1$-motive yield both equations
\eqref{eq:truncated-transfer} for the transported vectors in
$V_p(E)$ and classes in $\HdR(E/K)$, over the finitely generated
coefficient fields required here? Proposition~\ref{prop:MS-intersection}
specifies the distinction from the refined number-field
structures, and Proposition~\ref{prop:linear-truncation}
then supplies the passage to full vanishing.
\end{question}
\begin{question}
Can the unrestricted Drinfeld zero condition
$\mathcal E(z)=0$ of Proposition~\ref{prop:drinfeld-reformulation}
hold for an algebraic $z$ of degree at least three over $k$?
Corollary~\ref{cor:drinfeld-quadratic} excludes the lower degrees.
\end{question}

The full and normalized local relation ideals are given in
Corollary~\ref{cor:local-relation-ideal} and
Theorem~\ref{thm:graded-relation-ideal}.
Together with the filtration bound they describe what the
chosen local realization and its reductions detect.
Propositions~\ref{prop:motivic-A} and~\ref{prop:ordinary-arithmetic}
give the conditional consequences of the two period conjectures.
The complex and unrestricted Drinfeld questions, and the
arithmetic of non-CM Katz zeros over $\Q$, remain open.

\newcommand{\paperdoi}[1]{\href{https://doi.org/#1}{\textsc{doi}~\nolinkurl{#1}}}

\end{document}